\documentclass{article}
\usepackage{amsmath,amssymb,amsthm}
\usepackage[colorinlistoftodos,bordercolor=orange,backgroundcolor=orange!20,linecolor=orange,textsize=scriptsize]{todonotes}
\usepackage{fnpct}  
\newcommand{\R}{\mathbb{R}}

\newcommand{\E}[1]{\mathbf{E}\left[#1\right] }
\newcommand{\EE}[2]{\mathbf{E}_{#1}\left[#2\right] }

\newcommand{\norm}[1]{\lVert#1\rVert}

\newcommand{\Tr}[1]{\mbox{Tr}\left( #1\right)}

\DeclareMathOperator*{\argmin}{arg\,min}

\newtheorem{proposition}{Proposition}

\usepackage[margin=1.4in]{geometry}

\usepackage[backend=bibtex,citestyle=authoryear,maxcitenames=2,maxbibnames=5 ,url=false, doi=false]{biblatex}
\bibliography{../latex_files/full.bib}

\usepackage{url}            
\usepackage{booktabs}       
\usepackage{amsfonts}       
\usepackage{nicefrac}       
\usepackage{microtype}      
\usepackage{algorithm}
\usepackage[noend]{algpseudocode}
\usepackage{caption}
\usepackage{subcaption}
\graphicspath{{./images/}}

\title{Tracking the gradients using the Hessian:\\A new look at variance reducing stochastic methods}
\author{
  Robert M. Gower\\ LTCI, T\'el\'ecom-Paristech, Université Paris-Saclay\footnote{Robert M. Gower carried out this work while at INRIA in the Sierra team, funded by the Fondation de Sciences Math\'ematiques de Paris (FSMP) }\\\texttt{robert.gower@telecom-paristech.fr} \and
  Nicolas {Le Roux}\\Google Brain\\\texttt{nlr@google.com}  \and
  Francis Bach\\  INRIA, D\'epartement d'Informatique de l'ENS\\\texttt{francis.bach@inria.fr}
}

\begin{document}

\maketitle
\begin{abstract}
Our goal is to improve variance reducing stochastic methods through better control variates. We first propose a modification of SVRG which uses the Hessian to track gradients over time, rather than to recondition, increasing the correlation of the control variates and leading to faster theoretical convergence close to the optimum. We then propose accurate and computationally efficient approximations to the Hessian, both using a diagonal and a low-rank matrix. Finally, we demonstrate the effectiveness of our method on a wide range of problems.
\end{abstract}

\section{Introduction}
Many machine learning problems can be written as the minimization of a loss over a set of samples. Though this set can be infinite, it is in most cases finite, in which case the problem becomes the minimization of the average loss over $N$ samples, i.e.
\begin{align}
\theta^\ast &= \argmin_{\theta \in \R^d} \frac{1}{N}\sum_{i=1}^N f_i(\theta) = \argmin_{\theta\in \R^d} F(\theta),
\label{eq:argmin_theta}
\end{align}
where $f_i(\theta)$ is the loss incurred by parameters $\theta$ for the $i$-th sample. In this work, we make the assumption that each function $f_i$ is strongly convex.

A common way of solving Eq.~\eqref{eq:argmin_theta} is through a stochastic gradient method~\autocite{robbins1951stochastic}. Starting from $\theta_0$, the standard stochastic method samples $i$ uniformly at random in $\{1,\dots,N\}$ then computes $\theta_{t+1} = \theta_t - \gamma_t \frac{\partial f_i(\theta)}{\partial \theta}|_{\theta = \theta_t}$ where $\gamma_t$ can be a scalar or a full matrix. To unclutter notation, we will use $g$ to denote the gradient, i.e. $g_i(\theta_t) = \frac{\partial f_i(\theta)}{\partial \theta}|_{\theta = \theta_t}$.

Stochastic gradient updates are not convergent for constant stepsizes $\gamma_t$ since the variance of $g_i(\theta_t)$ does not converge to 0~\autocite{schmidt2014convergence}. Achieving convergence requires a decreasing schedule which impacts the convergence speed. By reducing the variance of these updates, one can hope to achieve convergence with larger or even constant stepsizes and thus obtain a faster rate. We now review the most popular stochastic variance reducing methods.

SAG~\autocite{leroux2012stochastic,schmidt2013minimizing} was the first stochastic method to achieve a linear convergence rate for strongly convex problems over a finite dataset. It required a storage linear in $N$ and used biased parameter updates. SDCA~\autocite{zhang2013linear,ShalevShwartz2013sdca} is a similar method solving a dual minimization problem for an important subclass of such problems. Compared to SAG which is adaptive, SDCA requires the knowledge of the strong convexity constant.

SVRG~\autocite{johnson2013accelerating} achieves a similar rate with a storage requirement independent of $N$, albeit at the expense of a constant factor increase in the computational cost. SVRG also requires an extra parameter besides the stepsize: the frequency at which the true gradient is recomputed. This is also the first work to provide the variance reduction interpretation. In this paper, we will mostly consider extensions of SVRG.

SAGA~\autocite{defazio2014saga} modifies SAG to get unbiased updates, leading to better convergence rates. Using a proximal algorithm, SAGA also works for non-smooth regularizers.

MISO~\autocite{mairal2013optimization} is a majorization-minimization algorithm achieving a linear convergence rate on strongly convex problems but it can also be used to solve nonconvex problems under some conditions.

All of these methods reduce variance by incorporating knowledge about gradients on all datapoints at every timestep rather than only using the gradient for the sampled datapoint. Where they differ is in which gradients they use and how they are incorporated.

This paper builds up on all these methods and makes the following contributions:
\begin{itemize}
\item We recall in Section~\ref{sec:control} that most of these methods can be seen as using control variates, an observation made by~\textcite{defazio2014saga}, then highlight that none of the control variates currently used track the gradients.
\item We propose in Section~\ref{sec:finite_second} a new control variate that tracks the gradient by building a linear model which is based on the Hessian matrix. This second-order information is used to further reduce the variance of the parameter updates, and not as a preconditioner.
\item Since calculating Hessians typically scales as $O(d^2)$, we propose several new approximations in Section~\ref{sec:approx} based on diagonal  or low-rank matrices. Such techniques are usually only dedicated to preconditioning.
\item We present in Section~\ref{sec:theory} a faster theoretical convergence rate than that of SVRG within a neighborhood of the solution.
\item Finally, in our experiments in Section~\ref{sec:experiments}, we show a significant and consistent improvement over SVRG.
\end{itemize}

Recently new work based on tracking the gradient using the Hessians of the $f_i$ functions has come to our attention~\autocite{WaiSNS17}. In this work the authors  
use incremental updates of the Hessian with the 2nd order Taylor expansion to maintain a rolling average estimate of the gradient. Our work also makes use of incremental Hessian updates together with the 2nd order Taylor expansion to
 estimate the gradient. Though differently from~\cite{WaiSNS17}  we focus on improving the SVRG method, as opposed to maintaining a rolling average.  Furthermore,~\cite{WaiSNS17}  use the full Hessian of the $f_i$ functions, which results in an iteration cost of $d^2$ while we focus on using diagonal and low rank Hessian approximations to keep the iteration cost linear in $d$.

\section{Control variates}
\label{sec:control}
We now recall that most variance reduced methods use control variates, a well-known technique in Monte-Carlo simulation (see, e.g.,~\textcite{rubinstein2016simulation}) designed to reduce the variance of the estimate of the expectation of a random variable. Let us start by describing this technique.

We try to estimate the unknown expectation $\bar{x}$ of a random variable $x$ and that we have access to another random variable, $z$, whose expectation $\bar{z}$ is known. Then the quantity $x_z = x - z + \bar{z}$ has expectation $\bar{x}$ and variance $V[x_z] = V[x] + V[z] - 2 \textrm{Cov}(x, z)$ where $V[]$ is the variance and $\textrm{Cov}[]$ the covariance. $V[x_z]$ is lower than $V[x]$ whenever $z$ is sufficiently positively correlated with $x$ and the variance reduction is larger when the control variate is more correlated with the random variable.

Using the true gradient $g(\theta_t) = \frac{1}{N} \sum_i g_i(\theta_t)$, batch methods achieve a linear convergence rate, albeit with an update computational cost linear in $N$. Stochastic methods estimate the true gradient through a single\footnote{Or a small number in the case of minibatch gradient methods.} sample $g_i(\theta_t)$, leading to a slower convergence rate but cheaper updates. The goal is thus to estimate the true gradient as accurately as possible while still maintaining update costs independent of $N$. As the true gradient is the expectation of the individual gradients, control variates can be applied.

More specifically, using a control variate to reduce the variance of the gradient estimation means replacing the quantity $g_i(\theta_t)$ with
\begin{align}
\tilde{g}_i(\theta_t) &= g_i(\theta_t) - z_i(\theta_t) + \frac{1}{N} \sum_j z_j(\theta_t) \; ,
\label{eq:control_variate}
\end{align}
transforming the stochastic gradient updates into
\begin{align}
\theta_{t+1} &= \theta_t - \gamma_t \Big(g_i(\theta_t) - z_i(\theta_t) + \frac{1}{N} \sum_j z_j(\theta_t) \Big)\; .
\label{eq:finite_update}
\end{align}
The new gradient estimate $\tilde{g}_i(\theta_t)$ is an unbiased estimator of the true gradient $g(\theta_t)$ with lower variance than $g_i(\theta_t)$ when $z_i(\theta_t)$ is positively correlated with $g_i(\theta_t)$.

The main difference between existing variance reducing gradient methods lies in the chosen $z_i(\theta_t)$. SAGA, for instance, chooses $z_i(\theta_t) = g_i(\bar{\theta}_i)$ the last gradient computed for sample $i$. Using such a control variate requires a storage linear in $N$. SVRG uses $z_i(\theta_t) = g_i(\bar{\theta}_t)$, the gradient computed at a fixed $\bar{\theta}_t$ for all samples. This allows for the computation of the control variate while only having to store $\bar{\theta}_t$. As it is the same for all datapoints, the storage requirement is independent of $N$.

One commonality to all existing methods is that each control variate remains constant until a new gradient has been computed for that datapoint. In particular, this means that control variates can easily become outdated if the gradients change quickly, leading to lower correlation and thus higher variance of the estimator. Since the variance of the estimator is directly linked to the convergence speed, improving correlation will lead to faster convergence. In the next section, we will show how to modify SVRG to make use of second-order approximations, maintaining highly correlated control variates longer than in SVRG.

\section{Tracking gradients with second-order control variates}
\label{sec:finite_second}
We now show how we can modify all stored gradients with a cost independent of $N$ by using second-order information about each individual function.

All algorithms of the previous section used for $z_i(\theta_t)$ the gradient computed for the same sample at a previous parameter $\bar{\theta}$. In particular, $z_i(\theta_t)$ does not change when the $i$th datapoint is not sampled. While this is efficient close to convergence when $\theta$ does not move much, correlation might be low early on, leading to slower convergence.

In particular, the convergence rate of SAGA is proportional to $\big(1 - \frac{\mu}{2 ( \mu N + L_{\max})}\big)$,  with $\mu$ the strong-convexity constant and $L_{\max}$ the maximum smoothness constant of the individual $f_i$'s.\footnote{For linear regression, we have $L_{\max} = \max_i \|x_i\|^2 = R^2$ the radius of the input data.} It is crucial to note that $L_{\max}$ is not the same as $L$, the smoothness constant of the full function $F$. More precisely, we have $L_{\max} /d \leq L \leq L_{\max}$. When all individual gradients are equal we have that $L = L_{\max}$ and SAGA is exactly $N$ times faster than full gradient descent. On the other hand, when the individual gradients are uniformly distributed over the whole space, we have $L = L_{\max} / d$ and the speed improvement of SAGA depends on $N / d$.

When $\mu/L_{\max}$ is much smaller than $1/N$, the poor condition number means that $\theta$ does not move much and SAGA should be close to batch gradient descent\footnote{The resulting algorithm would be $N$ times faster since an update only requires looking at a single datapoint compared to the full dataset for batch gradient descent.}. Indeed, the rate is close to $(1 - \mu/L_{\max})$ which is similar to $(1 - \mu/L)$, the convergence rate of batch gradient descent. However, the difference between $L$ and $L_{\max}$ means that, even in that regime, the rate of batch gradient descent can be much better. This is because batch gradient descent can make use of a stepsize of $1/L$, potentially much larger than the stepsize of $\frac{1}{2 (\mu N + L_{\max})}$ of SAGA.

When $\mu/L_{\max}$ is larger than $1/N$, the convergence rate of SAGA is close to $1 - 1/2N$, which is much smaller than the rate of convergence of batch gradient descent. This is because, in the well-conditioned case, $\theta$ changes quickly, leading to poorly correlated control variates and slower convergence relative to batch gradient descent.

Thus, finding control variates that have a greater correlation with the gradient would allow these algorithms to maintain a fast convergence for larger values of $N$ and potentially use larger stepsizes. Additionally, SVRG requires regular recomputation of $\bar{\theta}_t$ to maintain its convergence rate. Maintaining a high correlation for longer could decrease this required frequency of recomputation. 

Algorithm~\ref{alg:SVRG} describes the original SVRG method~\autocite{johnson2013accelerating}.
\begin{algorithm}[ht]
\begin{algorithmic}
\State \textbf{Parameter:} Functions $f_i$ for $i=1,\ldots, N$
\State Choose $\bar{\theta}_0 \in \R^d$ and stepsize $\gamma >0$
\For {$k = 0, \dots, K-1$}
    \State Calculate $g(\bar{\theta}_k)=\frac{1}{N}\sum_{j=1}^N g_j(\bar{\theta}_k)$, $\theta_0 = \bar{\theta}_k$
    \For {$t = 0, 1, 2, \dots, T-1$}
    \State $i \sim \mathcal{U}[1, N]$
	\State $\theta_{t+1} = \theta_t -\gamma \Big(g_i(\theta_t) - g_i(\bar{\theta}_k)  + g(\bar{\theta}_k)\Big) $
    \EndFor	
    \State $\bar{\theta}_{k+1} = \theta_T$
\EndFor
\State Output $\bar{\theta}_K$
\end{algorithmic}
\caption{SVRG}
\label{alg:SVRG}
\end{algorithm}
\begin{algorithm}[ht]
\begin{algorithmic}
\State \textbf{Parameter:} Functions $f_i$ for $i=1,\ldots, N$
\State Choose $\bar{\theta}_0 \in \R^d$ and stepsize $\gamma >0$
\For {$k = 0, \dots, K-1$}
    \State Calculate $g(\bar{\theta}_k)=\frac{1}{N}\sum_{j=1}^N g_j(\bar{\theta}_k)$, $\underline{H(\bar{\theta}_k)=\frac{1}{N}\sum_{j=1}^N H_j(\bar{\theta}_k)}$, $\theta_0 = \bar{\theta}_k$
    \For {$t = 0, 1, 2, \dots, T-1$}
    \State $i \sim \mathcal{U}[1, N]$
	\State $\theta_{t+1} = \theta_t -\gamma \Big(g_i(\theta_t) - g_i(\bar{\theta}_k)  \underline{- H_i(\bar{\theta}_k)(\theta_t - \bar{\theta}_k)}+ g(\bar{\theta}_k) + \underline{H(\bar{\theta}_k)(\theta_t - \bar{\theta}_k)}\Big) $
    \EndFor	
    \State $\bar{\theta}_{k+1} = \theta_T$
\EndFor
\State Output $\bar{\theta}_K$
\end{algorithmic}
\caption{SVRG2 (SVRG with tracking)}
\label{alg:SVRG2}
\end{algorithm}

We now modify Algorithm~\ref{alg:SVRG} by using control variates based on the first order Taylor expansion of the stochastic gradient. Rather than $z_j(\theta_t) = g_j(\bar{\theta})$, we will use:
\begin{align}
z_j(\theta_t) &= g_j(\bar{\theta}) + H_j(\bar{\theta})(\theta_t - \bar{\theta}),
\label{eq:z_t}
\end{align}
where $H_j(\bar{\theta})$ is the Hessian of $f_j$ taken at $\bar{\theta}$. The new control variate $z_j(\theta_t)$ now depends on $\theta_t$, which is not the case for the control variate used in SVRG.

Plugging Eq.~\eqref{eq:z_t} into Eq.~\eqref{eq:finite_update} yields the following update:
\begin{align}
\theta_{t+1} &= \theta_t -\gamma_t \Big(g_i(\theta_t) - g_i(\bar{\theta})
- H_i(\bar{\theta})(\theta_t - \bar{\theta}) 
+ \frac{1}{N}\sum_j g_j(\bar{\theta}) +\frac{1}{N}\sum_j H_j(\bar{\theta})(\theta_t - \bar{\theta}) \Big) \;.
\label{eq:osaga_update}
\end{align}
The resulting method, SVRG2, is presented in Algorithm~\ref{alg:SVRG2}. The only changes compared to Algorithm~\ref{alg:SVRG} are the computation of $H(\bar{\theta})$ and the terms involving the Hessian in the update (underlined). The theoretical analysis of Algorithm~\ref{alg:SVRG2} for generalized linear models is carried out in Section~\ref{sec:analysis}.

It is apparent that, if the $f_i$'s are exactly quadratic, then $g_i(\bar{\theta}) + H_i(\bar{\theta})(\theta_t - \bar{\theta}) = g_i(\theta_t)$ and the update is the same as batch gradient descent, albeit with a computational cost independent of $N$. The goal is thus to yield an update as close as possible to that of batch gradient descent without paying the price of recomputing all gradients every time.

Algorithm~\ref{alg:SVRG2} requires the computation of the full Hessian every time a new $\bar{\theta}$ is chosen and does two matrix-vector products for each parameter update. The cost of computing the full Hessian is $O(Nd^2)$ and the cost of each update is $O(d^2)$. Since $T$, the number of updates before recomputing the full gradient, is typically a multiple of $N$, the average cost per update is $O(d^2)$, which is $d$ times bigger than all first-order stochastic methods. As this cost is prohibitive for everything but low-dimensional problems, we need to resort to using approximations of the Hessian.

\section{Approximations of the Hessian}
\label{sec:approx}
Algorithm~\ref{alg:SVRG2} uses the following quantities:
\begin{itemize}
\item $H_i(\bar{\theta})$ needed for multiplying with the current~$\theta_t$;
\item $\sum_j H_j(\bar{\theta})$ needed for multiplying with the current~$\theta_t$;
\item $H_i(\bar{\theta})\theta_i$, $\sum_j g_j(\bar{\theta})$ and $\sum_j H_j(\bar{\theta})\bar{\theta}$ which are all vectors.
\end{itemize}
The first and second terms cannot be stored for large dimensional problems and, even in that case, the computational cost will be prohibitive. We thus need to resort to approximations $\hat{H}_j(\bar{\theta})$ of each Hessian $H_j(\bar{\theta})$ such that their sum $\hat{H}(\bar{\theta})$ is also easy to store.

This restriction makes us focus on linear projections of the Hessians, and two in particular: diagonal approximations and projections onto a low-rank subspace. Note that both of these classes of approximations include the special case $H_i = 0$, which is the original SVRG.

\subsection{Diagonal approximation}
We propose here a diagonal approximation of each Hessian leading to a storage cost of $O(d)$.

\textbf{Diagonal of the Hessian:}
By minimizing the Frobenius distance with the true Hessian, the diagonal of the Hessian tries to approximate the true Hessian in all directions of the space. This is in general very conservative as we care mostly about approximating the Hessian in the directions $\theta_t - \bar{\theta}$.

\textbf{Using the secant equation:}
Since the approximation is used in the term $\hat{H}_i(\bar{\theta})(\theta_t - \bar{\theta})$, we want $\hat{H}_i(\bar{\theta})(\theta_t - \bar{\theta})$ to approximate well $H_i(\bar{\theta})(\theta_t - \bar{\theta})$. One can use the secant equation and set
\begin{align*}
\hat{H}_i(\bar{\theta}) &= \frac{H_i(\bar{\theta})(\theta_t - \bar{\theta})}{\theta_t - \bar{\theta}} \approx \frac{g_i(\theta_t) - g_i(\bar{\theta})}{\theta_t - \bar{\theta}} \; ,
\end{align*}
where the division between two vectors is to be taken elementwise.
However, this leads to unstable updates as no guarantee is given outside of one direction.

\textbf{Robust secant equation:}
We can robustify the secant equation by minimizing the average squared-$\ell_2$ distance within a small ball around the previous direction. In other words, we seek
\begin{align*}
\hat{H} &= \argmin_Z \int_{\Delta} \left\|Z(\theta_t - \bar{\theta} + \Delta) - H_i(\bar{\theta})(\theta_t - \bar{\theta} + \Delta)\right\|^2 p(\Delta) \; d\Delta.
\end{align*}
Assuming $\Delta \sim \mathcal{N}(0, \sigma^2 I)$, we get
\begin{align}
\hat{H} &= \frac{(\theta_t - \bar{\theta})\odot H_i(\bar{\theta})(\theta_t - \bar{\theta}) + \sigma^2 \textrm{diag}(H_i(\bar{\theta}))}{(\theta_t - \bar{\theta}) \odot (\theta_t - \bar{\theta}) + \sigma^2}\approx\frac{(\theta_t - \bar{\theta}) \odot (g_i(\theta_t) - g_i(\bar{\theta})) + \sigma^2 \textrm{diag}(H_i(\bar{\theta}))}{(\theta_t - \bar{\theta}) \odot (\theta_t - \bar{\theta}) + \sigma^2}.\label{eq:robsec}
\end{align}
This is an interpolation between the two previous approximations. In particular, $\sigma^2 = +\infty$ recovers the diagonal of the Hessian while $\sigma^2 = 0$ recovers the secant equation. For a fixed $\sigma^2$, larger parameters updates at the beginning of optimization will make this approximation close to the secant equation while, close to convergence, this approximation will converge towards the diagonal of the Hessian due to smaller parameter updates.

\subsection{Low-rank approximation: Curvature Matching}
Another possibility is to build a low rank approximation of the Hessian by using a low dimensional embedding. For brevity let $H_i := H_i(\bar{\theta})$ and $H :=H(\bar{\theta}).$

Consider a matrix $S \in \R^{d \times k}$ where $k \ll d, N$. The matrix $S$ may be random. Calculating the \emph{embedded} Hessian $S^\top H_i S$ gives us a low rank approximation of $H_i$  using the following model:
\begin{align}
\hat{H}_i & =  \argmin_{X\in \R^{d \times d}} \|X \|_{F(H)}^2 \nonumber \\
&\quad \mbox{subject to} \quad S^\top XS = S^\top H_iS,\label{eq:curvmatch}
\end{align}
where $\|\hat{H} \|_{F(H)}^2 := \Tr{\hat{H}^\top H\hat{H} H}$ is the weighted Frobenius norm. The solution to the above is a rank $k$ matrix given by
\begin{equation}\label{eq:embedleftright}
 \hat{H}_i = HS(S^\top H S)^{\dagger}S^\top H_iS(S^\top H S)^{\dagger}S^\top H\in \R^{d\times d},
 \end{equation}
where $M^\dagger$ is the pseudoinverse of $M$. We refer to methods based on the Hessian approximation presented in Eq.~\eqref{eq:embedleftright} as the \emph{Curvature Matching} methods since the approximation has the same curvature as the true Hessian over the subspace spanned by $S$.

Our use of a Hessian weighted Frobenius norm in~\eqref{eq:curvmatch} is largely inspired on quasi-Newton methods. Indeed, during the development of quasi-Newton methods it was observed that using the Hessian as a weighting matrix in the Frobenius norm gave rise to a more efficient method~\autocite{Goldfarb1970}, namely the BFGS method. Our findings mirror that of the quasi-Newton literature in that our empirical tests showed that using the Hessian weighted Frobenius norm is more efficient than using the standard Frobenius norm.

Note that, to compute Eq.~\eqref{eq:embedleftright}, we only need to store the $d\times k$ matrix $H S$ and calculate $S^\top H_iS$. The advantages of this approach are numerous and include: (a) Reduced memory storage, (b) Hessian-vector products $\hat{H}_{i}v$ can be computed efficiently at a cost of $O(k(  2d +3k))$ which is linear in $d$, (c) the approximation $\hat{H}_{i}$ is the result of applying a linear operator on $H_i$ and consequently the expectation of $\hat{H}_{i}$ can be computed efficiently.

This last item is of particular importance for maintaining an unbiased estimator. Indeed, taking the expectation  over~\eqref{eq:embedleftright}, we have that
\begin{eqnarray}
 \EE{i}{ \hat{H}_i} &  \overset{\eqref{eq:embedleftright}}{=}& HS(S^\top H S)^{\dagger}S^\top \EE{i}{H_i}S(S^\top H S)^{\dagger}S^\top H \nonumber \\
 &=& HS(S^\top H S)^{\dagger}S^\top H S(S^\top H S)^{\dagger}S^\top H \nonumber  \\
 &= &  HS(S^\top H S)^{\dagger}S^\top H,\label{eq:EJHjcurv}
\end{eqnarray}
where in the last step we used that $M^\dagger M M^\dagger =M^\dagger$ for all matrices $M$. From Eq.~\eqref{eq:EJHjcurv} we see that the expectation of $\hat{H}_i$ is equal to the resulting approximation were we to use $H$ in place of $H_i$ in Eq.~\eqref{eq:embedleftright}. Furthermore~\eqref{eq:EJHjcurv} is a low rank matrix which can be stored and computed efficiently since we only need to calculate the action of the true Hessian $HS$.\footnote{With automatic differentiation, it costs only $k$ times the cost of evaluating the function $F(\theta)$ to compute $HS.$}

Specifically, in the outer iterations, we can compute and store the $d \times k$ Hessian \emph{action} matrix $A =\frac{1}{N}\sum_{j=1}^N H_j(\bar{\theta})S$ and the $k\times k$ Hessian \emph{curvature} matrix $C =(S^\top\frac{1}{N}\sum_{j=1}^N H_j(\bar{\theta})S)^{\dagger/2}$  then use $A$ and $C$ to perform the necessary matrix vector products. See Algorithm~\ref{alg:SVRGemb} for our implementation. Again we underline the computations that need to be added to the SVRG algorithm to produce our new algorithm.

\begin{algorithm}[ht]
\begin{algorithmic}
\State \textbf{Parameter:} Functions $f_i$ for $i=1,\ldots, N$
\State Choose $\bar{\theta}_0 \in \R^d$ and stepsize $\gamma >0$
\For {$k = 0, \dots, K-1$}
    \State Calculate $g(\bar{\theta}_k)=\frac{1}{N}\sum_{j=1}^N g_j(\bar{\theta}_k)$, $\theta_0 = \bar{\theta}_k$
    \State Calculate  $\underline{A = \frac{1}{N}\sum_{j=1}^N H_j(\bar{\theta}) S}$, $\underline{C = (S^\top A)^{\dagger/2}}.$
    \State Generate $S \in \R^{d \times k}$, calculate $\underline{\bar{S} = SC}$.
    \State Normalize the Hessian action $\underline{\bar{A} = AC}.$
    \For {$t = 0, 1, 2, \dots, T-1$}
    \State $i \sim \mathcal{U}[1, N]$
    \State $d_t = g_i(\theta_t) - g_i(\bar{\theta}_k) + g(\bar{\theta}_k)-\underline{ \bar{A}\bar{S}^\top H_i(\bar{\theta}_k)\bar{S}\bar{A}^\top  (\theta_t - \bar{\theta}_k)}+\underline{\bar{A}\bar{A}^\top (\theta_t - \bar{\theta}_k)}$
	\State $\theta_{t+1} = \theta_t -\gamma d_t$
    \EndFor	
    \State $\bar{\theta}_{k+1} = \theta_T$
    \State Output $\bar{\theta}_{K}$	
\EndFor
\end{algorithmic}
\caption{CM: Curvature Matching}
\label{alg:SVRGemb}
\end{algorithm}
%

While $S$ can be any matrix, including a random one drawn from some distribution, 
we found the method performed best when the construction of $S$ was based on the step directions $d_t$ for $t = 0 ,\ldots, T-1$ taken from the inner loop. 
The reasoning behind this choice is that, since ultimately we only require the action of the Hessian approximation along the directions $\theta_t - \bar{\theta}$,  the matrix $S$ should be chosen in such a way that our Hessian approximation is accurate along these directions. Thus we will construct the embedding matrix $S$ by again taking inspiration from the BFGS~\autocite{Goldfarb1970} and L-BFGS~\autocite{nocedal1980} methods and column concatenate previous step directions. Specifically, since there are many step directions and sequential step directions are often correlated, we arranged  the $T$ step directions into $k$ sets. We used the average over each of these $k$ sets as a column of $S$. That is
\begin{equation}\label{eq:prevS}
S = [\bar{d}_0, \ldots, \bar{d}_{k-1}] \quad \mbox{where} \quad \bar{d}_i =\frac{k}{T} \sum_{j= \frac{T}{k}i}^{\frac{T}{k}(i+1)-1} d_j,\end{equation}
where we assume, for simplicity, that $k$  divides $T.$

\subsubsection{Computational complexity}
Algorithm~\ref{alg:SVRGemb} can be applied to large-scale problems due to its reduced complexity. Indeed, using automatic differentiation~\autocite{Walther2008,Christianson:1992} the cost of computing each one of the matrix-matrix products $H_j(\bar{\theta}) S$ is $O(k \cdot \mbox{eval}(f_j))$, where $\mbox{eval}(f_j)$ is the cost of evaluating the function $f_j$\footnote{In fact, these matrix-matrix products are computed using directional derivatives, and it is also easy enough to handcode these directional derivatives to achieve this complexity.}. Furthermore, evaluating a single gradient $g_j$ also costs $O(\mbox{eval}(f_j)).$ The matrix product $S^\top A$ costs $O(k^2 d)$ to compute, while the products $AC$ and $SC$ cost $dk^2.$ 
 Computing the square and pseudoinverse in $C = (S^\top A)^{\dagger/2}$ costs $O(k^3)$, leading to a total for the outer iteration of
  \[O(k \cdot\sum_{j=1}^n \mbox{eval}(f_j) + k^2d + k^3) = O(k \cdot\mbox{eval}(f) + k^2d + k^3).\]
The main cost in the inner iteration is computing  $H_i(\bar{\theta}_k) S$ and a handful of matrix-vector products where the matrices have dimensions $k \times d$ or $d \times k$. The total cost of each inner iteration is $O(k\cdot \mbox{eval}(f_i)+ dk^2)$. Overall, the cost of computing an iteration (inner or outer) of Algorithm~\ref{alg:SVRGemb} is linear in $d$ and comes at a $k$-fold increase of the cost of computing an SVRG iteration. In our experiments, we used $k=10$.

\subsection{Low-rank approximation: Action Matching}

In the previous section we developed a low rank update based on curvature matching, though ultimately we required the action of the true Hessian $HS$ to compute the resulting Hessian approximation~\eqref{eq:embedleftright}.
 With the action of the Hessian  we can also build an approximation of the Hessian using an \emph{action matching} model:
\begin{align}
\hat{H}_i & =  \argmin_{X\in \R^{d \times d}} \,\norm{X}_{F(H)}^2\nonumber\\
 &\quad \mbox{subject to} \quad XS =  H_iS, \quad X = X^\top.\label{eq:actionmatch}
\end{align}
In other words, $\hat{H}_i$ is the symmetric matrix with the smallest norm which matches the action of the true Hessian.
The solution to the above is a rank $2k$ matrix given by
\begin{align}
\hat{H}_i &= HS( S^T  H S)^{-1}S^\top H_i\left(I-S( S^T  H S)^{-1}S^\top H\right) + H_iS( S^T  H S)^{-1}S^\top H.\label{eq:DFPinv}
\end{align}
Differently from the curvature matching model~\eqref{eq:curvmatch}, here we have had to explicitly enforce the symmetry of the Hessian approximation as a constraint~\eqref{eq:actionmatch} as the solution of the unconstrained model might not be symmetric.

This approximation is a generalization of the Davidson-Fletcher-Powell (\cite{Fletcher1960,Davidon1968}) update. This update, and a regularized inverse of this update, have already been used in an optimization context (\cite{Gower2014c}) but only as a preconditioner.

Again, we have that the approximation~\eqref{eq:DFPinv} is the result of applying a linear function to $H_i$. Thus, taking the expectation, we have that
\begin{align*}
\EE{i}{ \hat{H}_i}
&= HS( S^T  H S)^{-1}S^\top \EE{i}{H_i} +\EE{i}{ H_i}S( S^T  H S)^{-1}S^\top H\\
& \quad - HS( S^T  H S)^{-1}S^\top  \EE{i}{H_i} S( S^T  H S)^{-1}S^\top H \\
&=HS( S^T  H S)^{-1}S^\top H,
\end{align*}
which is the same as~\eqref{eq:EJHjcurv} and thus can be efficiently computed.

 The pseudo-code of our implementation of the method based on~\eqref{eq:actionmatch} is in Algorithm~\ref{alg:actionmatch}. The total complexity of the Algorithm~\ref{alg:actionmatch} is the same as the complexity of Algorithm~\ref{alg:SVRGemb}, though the hidden constant in the complexity is a bit larger as there is an additional matrix-vector product.

\begin{algorithm}[ht]
\begin{algorithmic}
\State \textbf{Parameter:} Functions $f_i$ for $i=1,\ldots, N$
\State Choose $\bar{\theta}_0 \in \R^d$ and stepsize $\gamma >0$
\For {$k = 0, \dots, K-1$}
    \State Calculate $g(\bar{\theta}_k)=\frac{1}{N}\sum_{j=1}^N g_j(\bar{\theta}_k)$, $\theta_0 = \bar{\theta}_k$
    \State Calculate and store $\underline{A = \frac{1}{N} \sum_{j=1}^N H_j(\bar{\theta}) S}$ and $\underline{C = (S^\top A)^{\dagger/2}}.$
    \State Generate $S \in \R^{d \times k}$, calculate and store $\underline{\bar{S} = SC}$.
    \State Normalize the Hessian action $\underline{\bar{A} = AC}.$
    \For {$t = 0, 1, 2, \dots, T-1$}
    \State $i \sim \mathcal{U}[1, N]$
	\State $\theta_{t+1} = \theta_t -\gamma \Big(g_i(\theta_t) - g_i(\bar{\theta}_k) + g(\bar{\theta}_k)$ \\ $\qquad \qquad \qquad 
	-\underline{\left(\bar{A}\bar{S}^\top H_i\left(I-\bar{S}\bar{A}^\top\right)
  + H_i\bar{S} \bar{A}^\top\right)
	 (\theta_t - \bar{\theta}_k)}
	 +\underline{\bar{A}\bar{A}^\top (\theta_t - \bar{\theta}_k)}\Big)$
    \EndFor	
    \State $\bar{\theta}_{k+1} = \theta_T$
    \State Output $\bar{\theta}_{K}$	
\EndFor
\end{algorithmic}
\caption{AM: Action Matching}
\label{alg:actionmatch}
\end{algorithm}
A similar approximation to~\eqref{eq:actionmatch} has been developed before~\autocite{GowerGold2016,Gower2016} but focused on obtaining an approximation to the inverse Hessian to then be applied as a quasi-Newton update.


\subsection{Comparing approximations}
Both the diagonal and the low-rank approximations have a free parameter: $\sigma^2$ for the diagonal approximation and the rank $k$ for the low-rank approximation. These two parameters make  different tradeoffs. $\sigma^2$ controls the robustness of the approximation: a larger value is more stable but likely to decrease the convergence rate. $k$ controls the complexity of the approximation: a larger value leads to a better convergence rate in terms of updates but with a larger computational cost per update. Though this was not explored in this work, it seems natural to combine these two approximations.

\section{Theoretical Analysis and Discussion}
\label{sec:theory}
\textbf{Convergence rates:}
We now provide a convergence analysis for generalized linear models when an approximation of the  Hessian is used.  The full proofs of the results below may be found in the Appendix.
\label{sec:analysis}
We make the following assumptions:
\begin{itemize}
\item We consider a compact set $\Theta$ included in the ball of center $\theta_\ast$ and radius $D$, where $\theta_\ast$ is the minimizer of $F(\theta) = \frac{1}{N} \sum_{i=1}^N f_i(\theta)$.
\item $F$ is $\mu$-strongly convex and $L$-smooth on $\Theta$.
\item Each $f_i$ is of the form $f_i(\theta) = \varphi_i(x_i^\top \theta)$ with $\| x_i\|^2 \leqslant R^2$, and $\varphi_i''(x_i^\top \theta) \in [ \alpha, 1]$, $\varphi_i'''(x_i^\top \theta) \leqslant \beta$ for all $\theta \in \Theta$. In particular, this means that we can use $L_{\max} = R^2$.
\item We consider the control variate
$z_i(\theta) = f_i'(\bar\theta) + H_i ( \theta - \bar\theta)$, where $H_i$ is an approximation of the Hessian $ f''_i(\bar\theta)$,
with a relative error so that $\frac{1}{N} \sum_{i=1}^N ( f''_i(\bar\theta) - H_i)^2 \preccurlyeq R^2 \eta \frac{1}{N} \sum_{i=1}^n f''_i(\bar\theta) $. Note that $\eta=1$ if we take $H_i = 0$ (plain SVRG) and $\eta=0$ if we use the exact Hessian.
\end{itemize}
We obtained the following proposition which shows that we get the actual condition number $L/\mu$ rather than $L_{\max} / \mu$ (which is the one for SVRG), if we are close enough to the global optimum and the Hessians are sufficiently well approximated:
\begin{proposition}
Assume $  D^2 \leqslant \frac{L \alpha}{8 \beta^2 R^2}$,  $\gamma = 1/ (4L)$
with
$  \eta \leqslant  \frac{L \alpha}{8 \beta^2 } $. Then, after a single epoch of SVRG2 with $T \geqslant \frac{16 L}{\mu}$, the error is reduced by a factor of $3/4$.
\end{proposition}
The result above implies that, in order to reach a precistion $\varepsilon$, we need $K = O( \log \frac{1}{\varepsilon})$ epochs of SVRG2, with an overall number of accesses to gradient and Hessian oracles less than $ O ( N + \frac{L}{\mu}  \log \frac{1}{\varepsilon} )$, once we are close enough to the optimum (but note that the bound on $D$ does not depend on $\mu$). In the Appendix, we also show that our algorithm is robust even far away from optimum with a step-size $\gamma \propto 1/R^2$.

\textbf{Difference with reconditioning:}
To our knowledge, this is the first use of the Hessian to track gradients over time as it is usually used for reconditioning. While using the Hessian to track the gradients does not remove the dependency on the condition number, it has several advantages. First, since there is no matrix inversion, it is much more robust to approximations of the true Hessian, as shown in the theoretical analysis above. Second, while classical second-order methods only enjoy superlinear convergence close to the optimum, tracking gradients improves convergence in a larger ball around the optimum.

This does not preclude the use of reconditioning and we believe the Hessian could be used for both purposes at the same time. It remains to be seen, however, if the best Hessian approximation is the same for both uses.

\textbf{Difference to non-uniform sampling}
Non-uniform sampling aims to close the gap between $L_{max}$ and $\tfrac{1}{n}\sum_i L_i$~(\textcite{Kern:2016,SchmidtBADCS15}). Our method closes the gap between $L_{max}$ and $L$. These two quantities, $\tfrac{1}{n}\sum_i L_i$ and $L$, can be very different. For instance, for a linear least squares objective $f(w) = 1/(2n)|| Xw-b||_2^2$ we have that $\tfrac{1}{n}\sum_i L_i = (1/n)\Tr{XX^\top}$ while $L = \tfrac{1}{n}\lambda_{\max}(XX^\top)$, thus $L$ can be upto $n$ times smaller than  $\tfrac{1}{n}\sum_i L_i$.


\section{Experiments}
\label{sec:experiments}

\begin{figure}
	\centering
	\begin{subfigure}[t]{0.32\textwidth}
		\centering
		\includegraphics[width =  \textwidth]{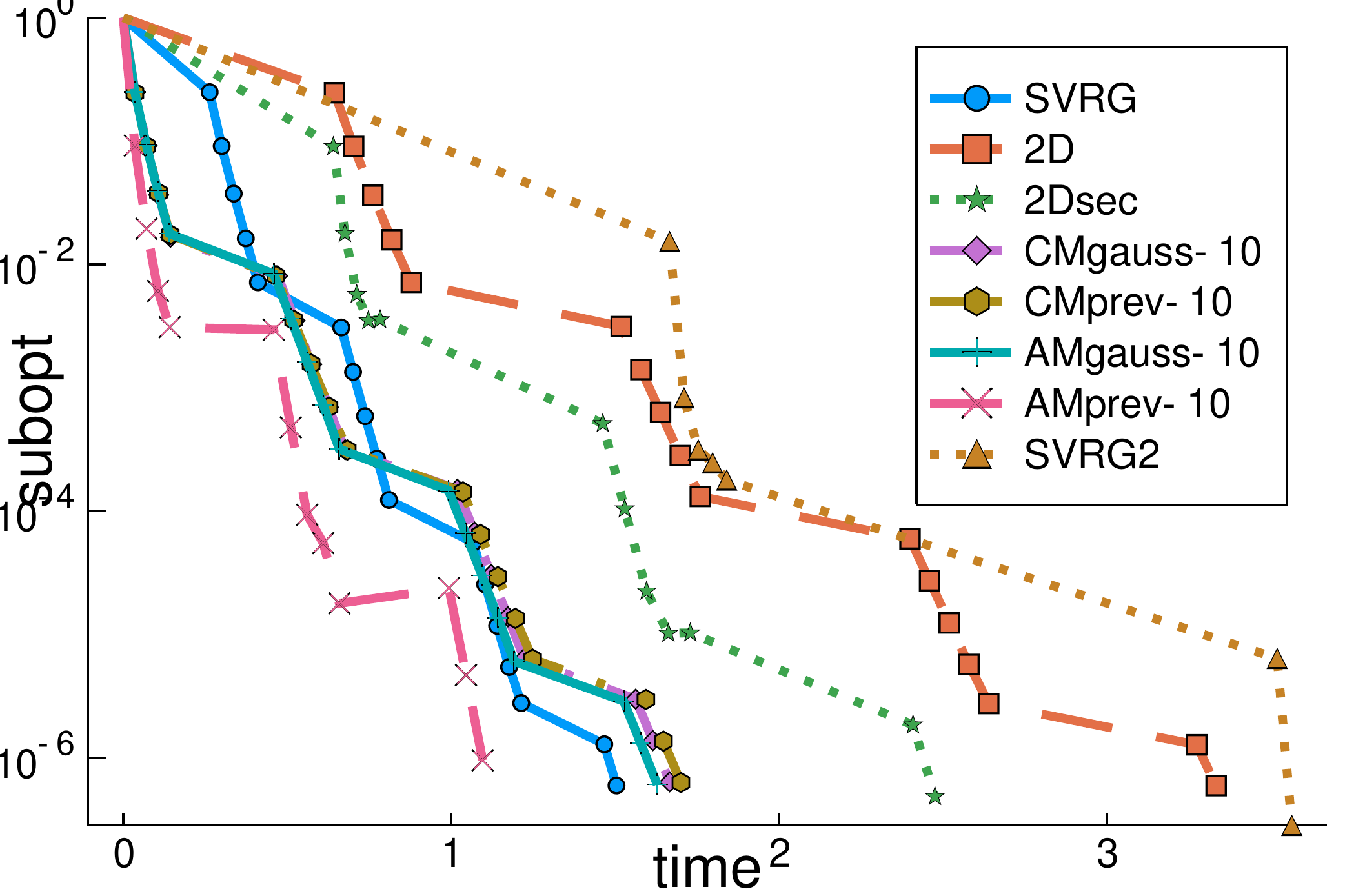} 
		\includegraphics[width =  \textwidth]{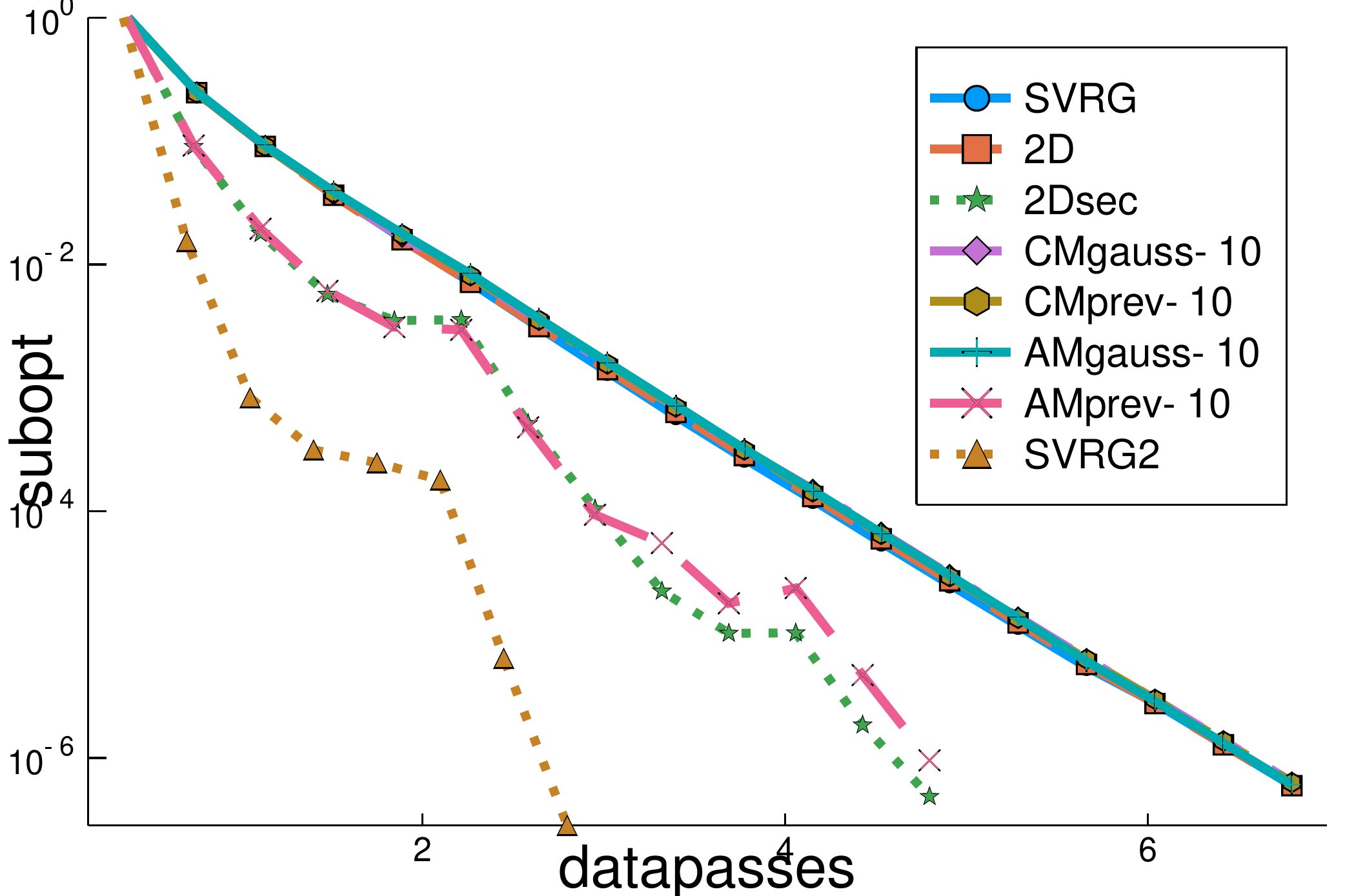}
		\caption{\texttt{a9a}}
	\end{subfigure}%
	\hspace{0.005\textwidth}
	\begin{subfigure}[t]{0.32\textwidth}
		\centering
		\includegraphics[width =  \textwidth]{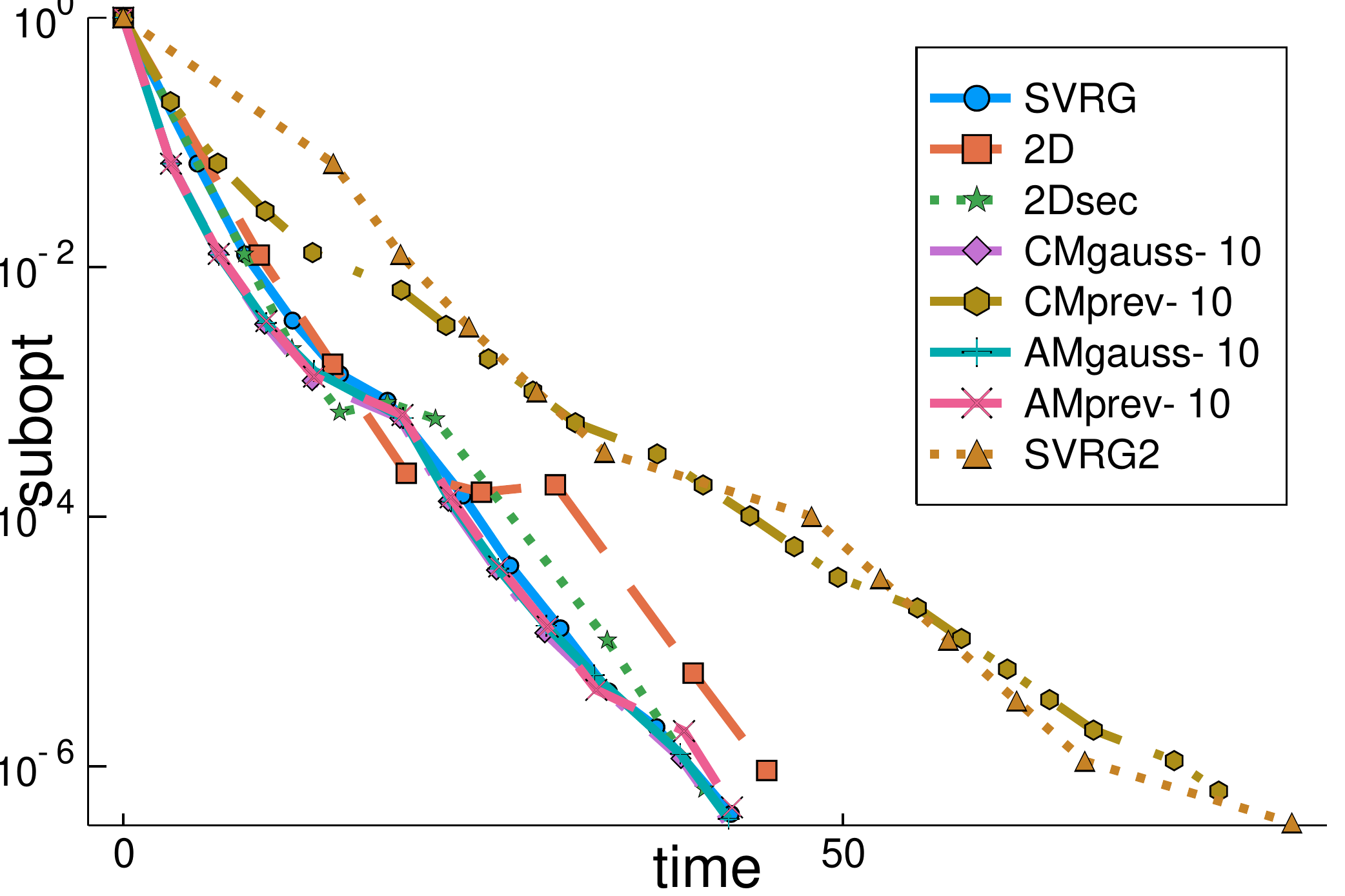} 
		\includegraphics[width =  \textwidth]{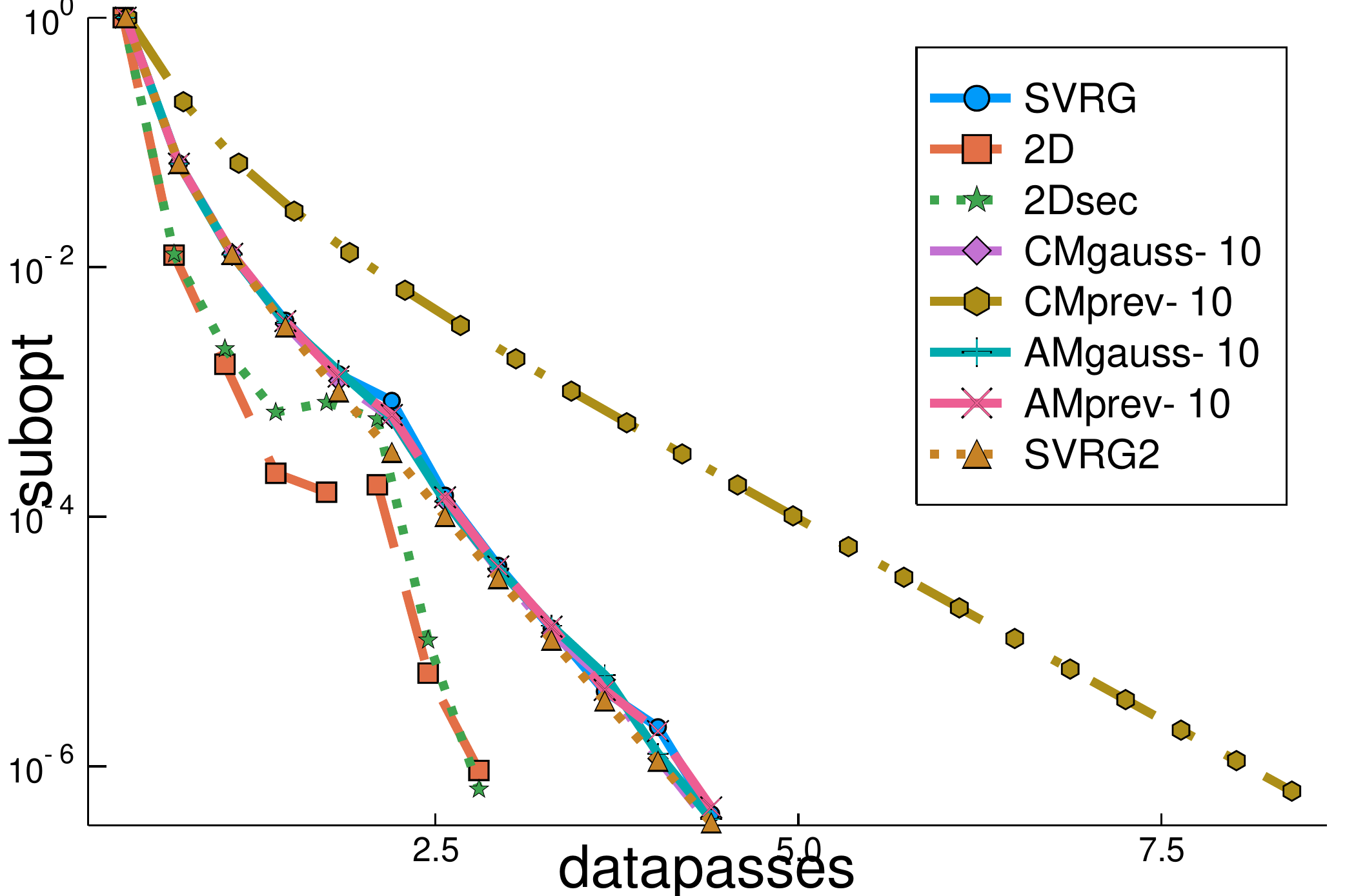}
		\caption{\texttt{covtype}}
	\end{subfigure}
	\hspace{0.005\textwidth}
	\begin{subfigure}[t]{0.32\textwidth}
		\centering
		\includegraphics[width =  \textwidth]{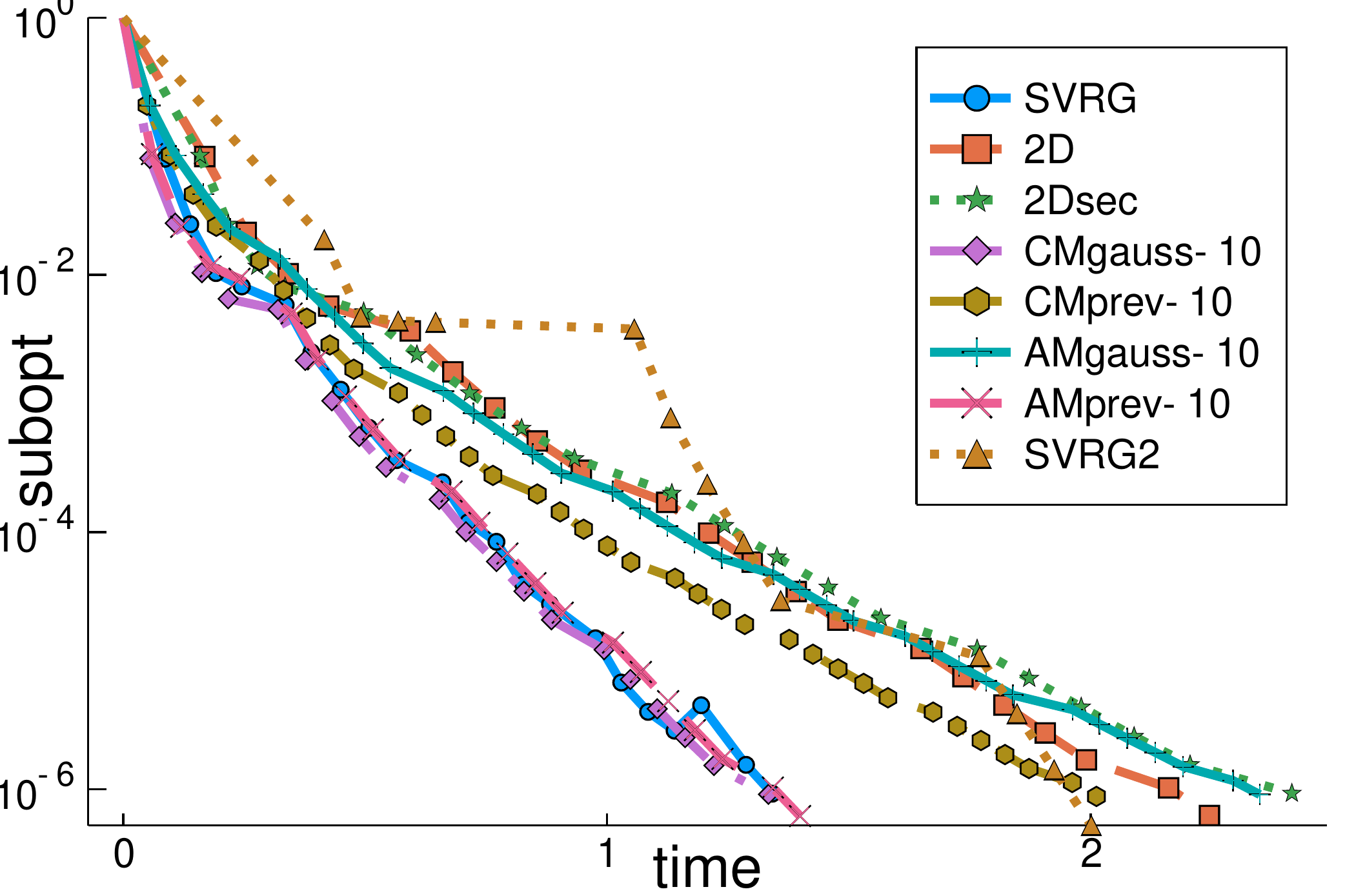} 
		\includegraphics[width =  \textwidth]{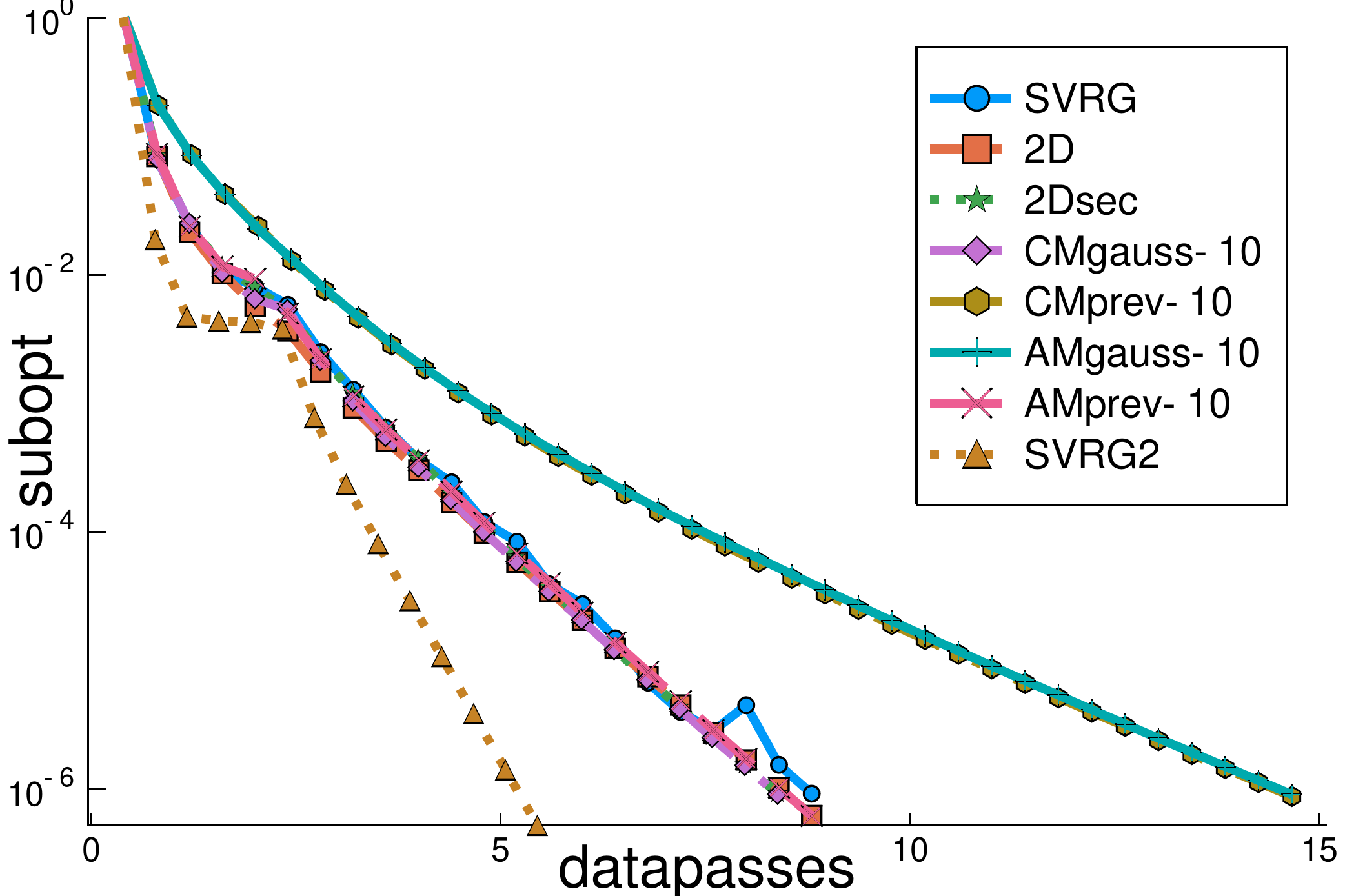}
		\caption{\texttt{mushrooms}}
	\end{subfigure}\\
	\begin{subfigure}[t]{0.32\textwidth}
		\centering
		\includegraphics[width =  \textwidth]{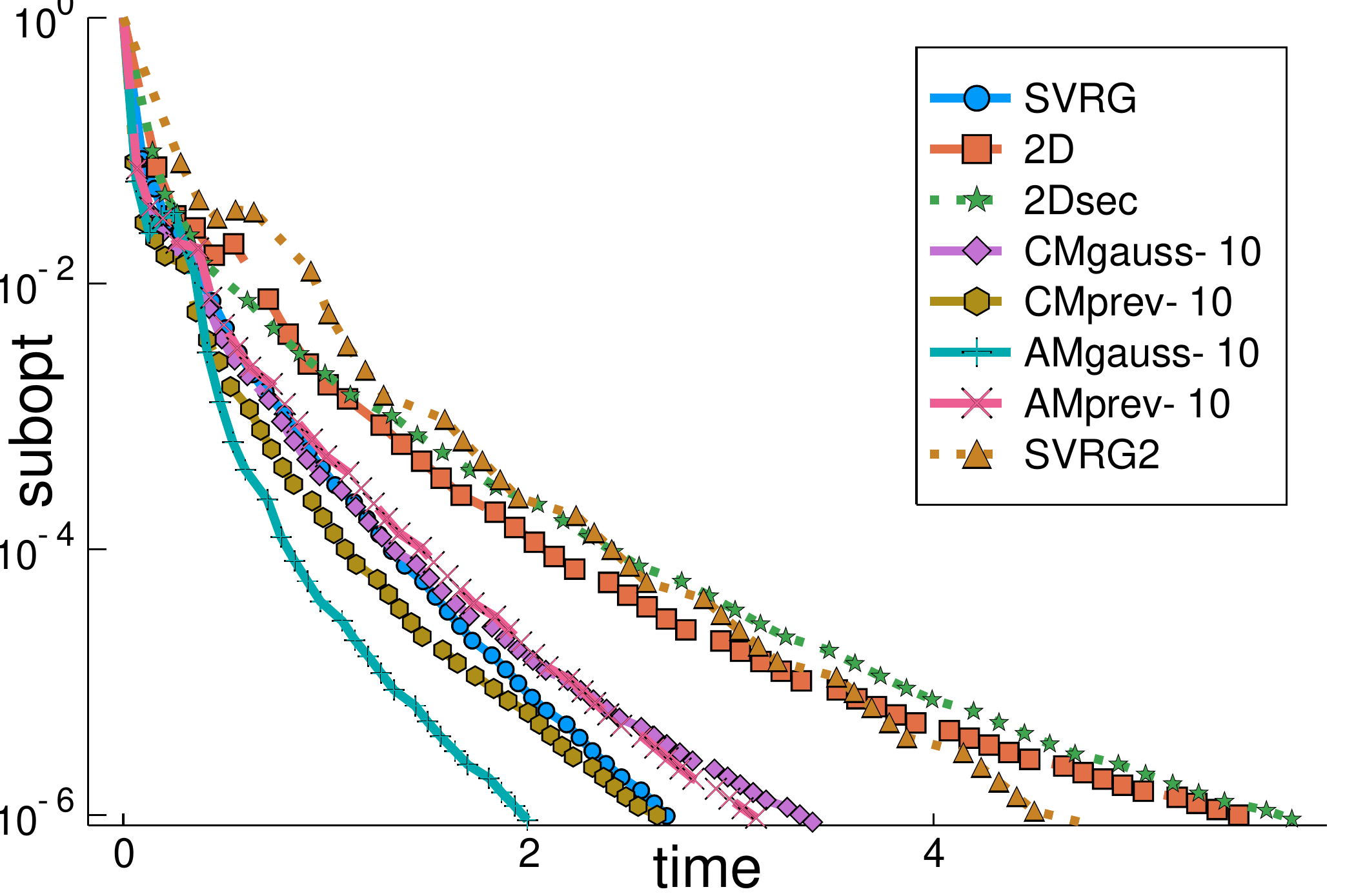} 
		\includegraphics[width =  \textwidth]{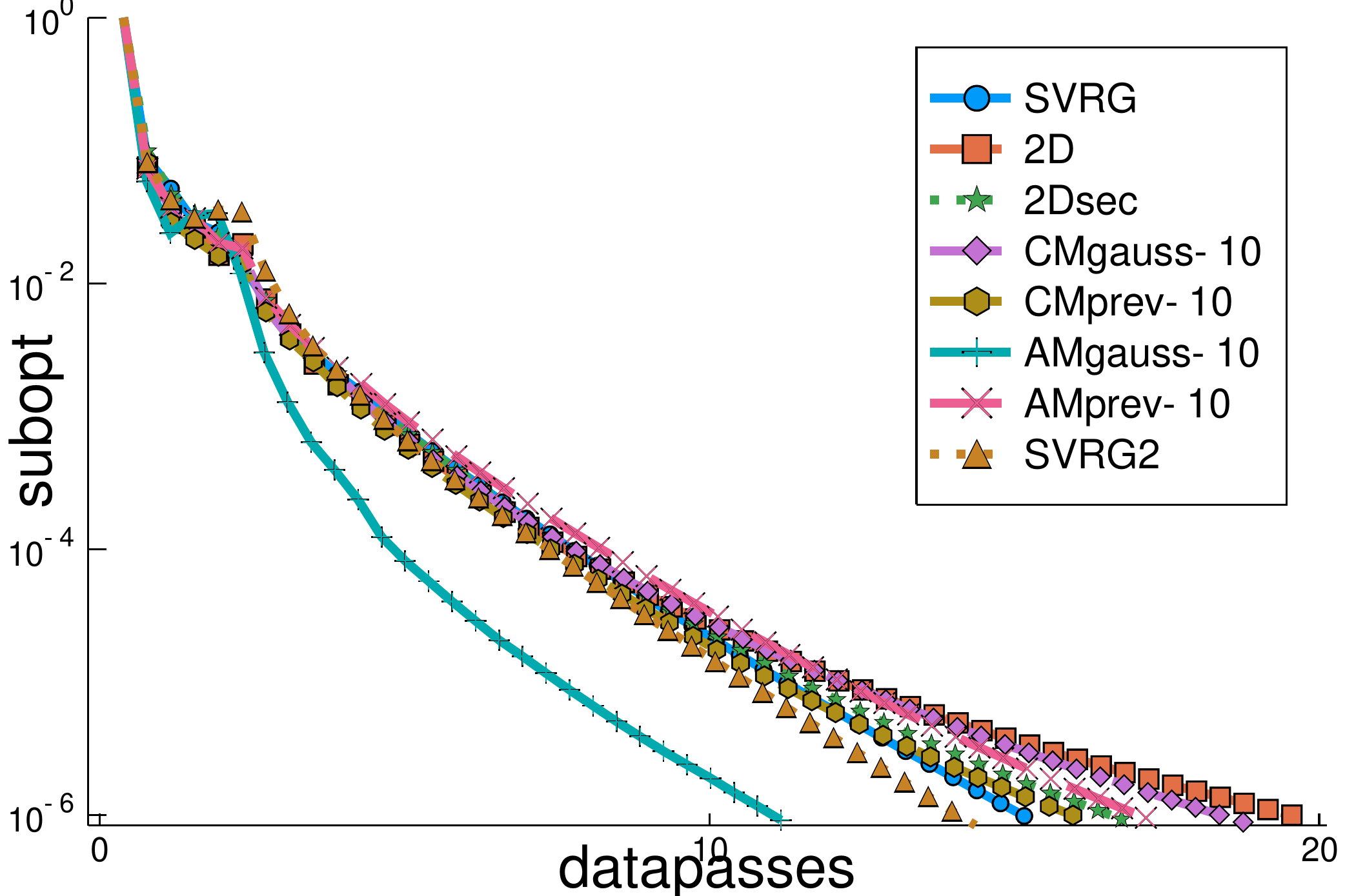}
		\caption{\texttt{phishing}}
	\end{subfigure}
	\hspace{0.005\textwidth}
	\begin{subfigure}[t]{0.32\textwidth}
		\centering
		\includegraphics[width =  \textwidth]{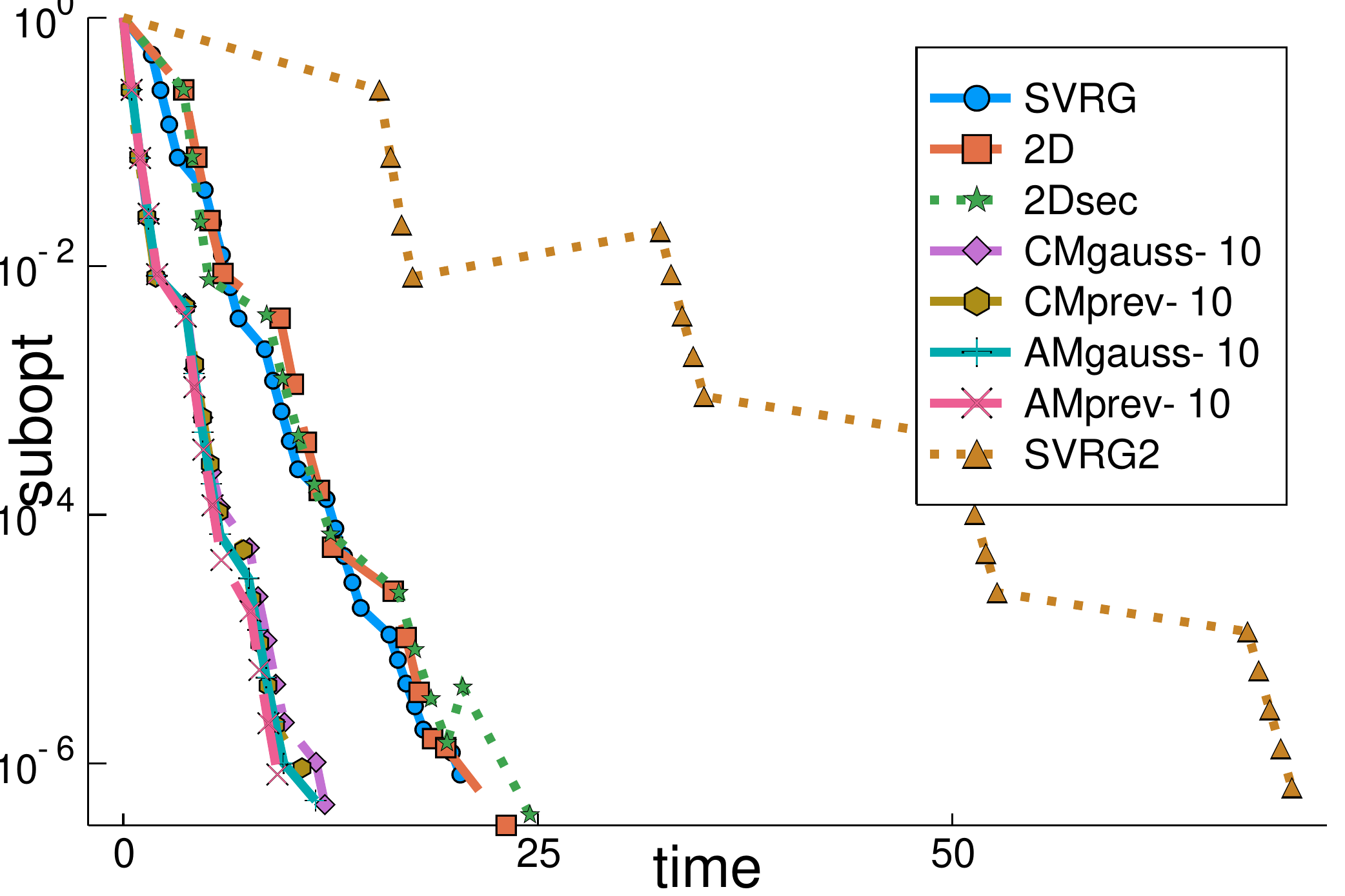} 
		\includegraphics[width =  \textwidth]{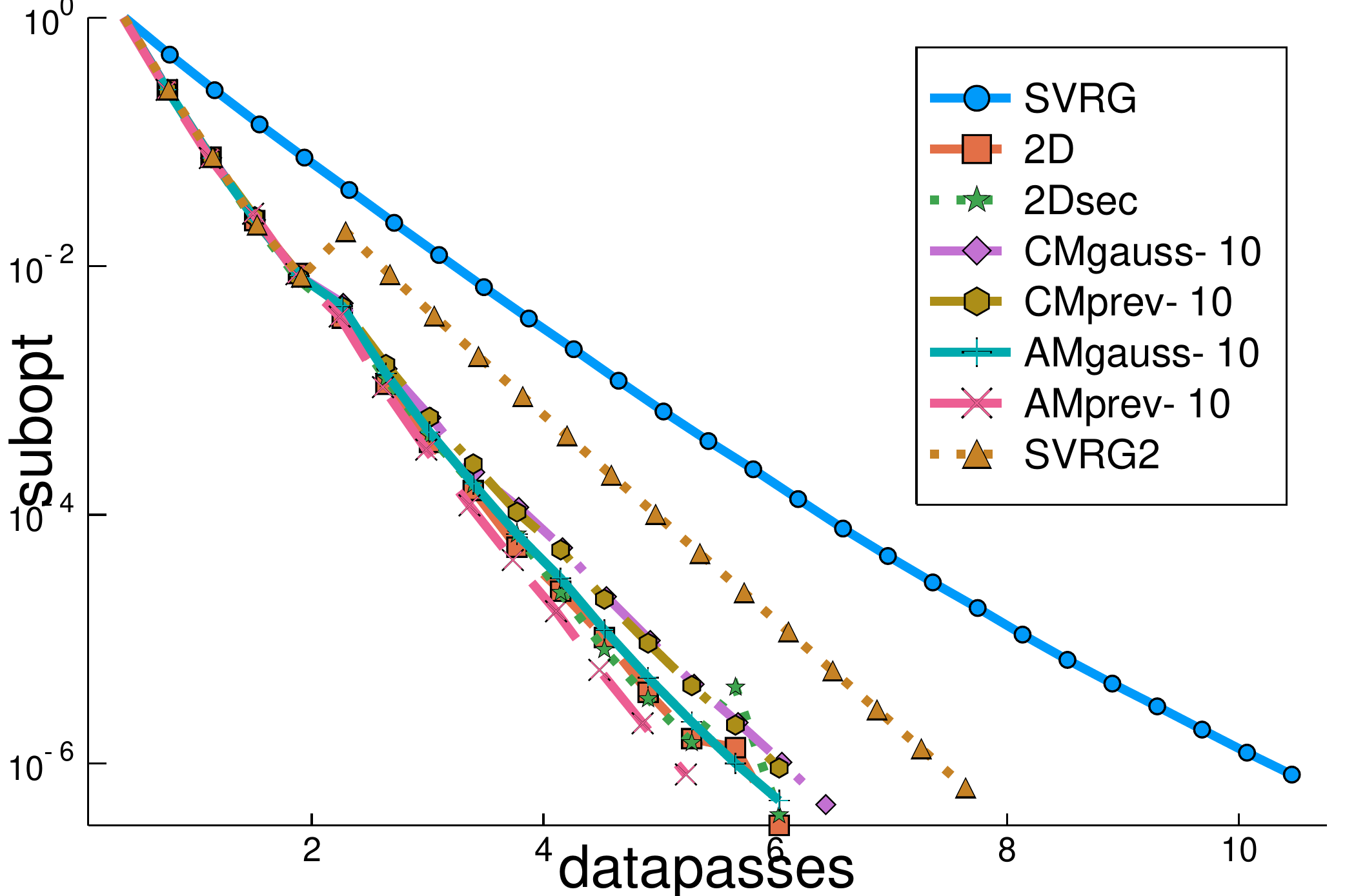}
		\caption{\texttt{w8a}}
	\end{subfigure}
	\hspace{0.005\textwidth}
	\begin{subfigure}[t]{0.32\textwidth}
		\centering
		\includegraphics[width =  \textwidth]{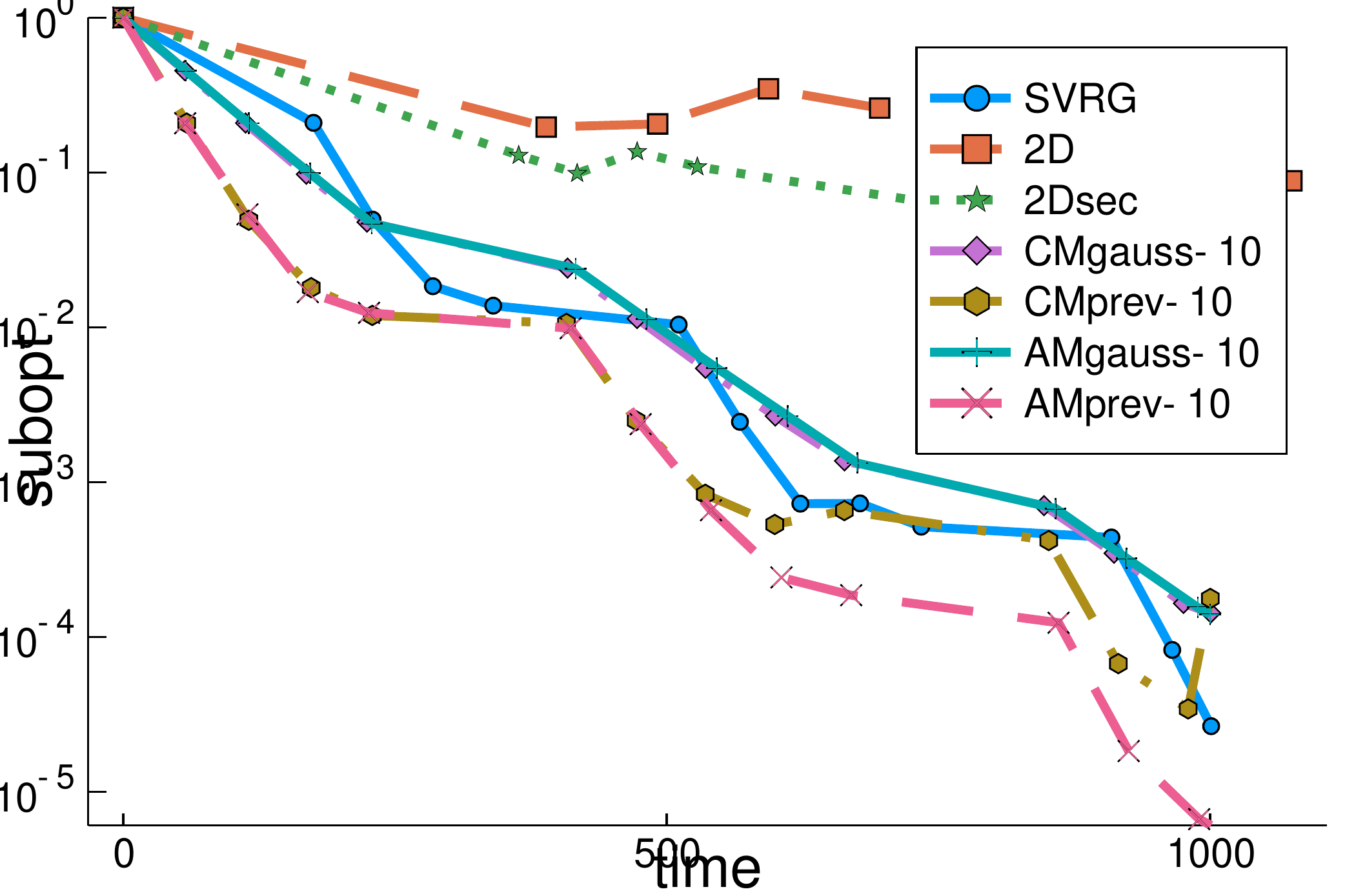} 
		\includegraphics[width =  \textwidth]{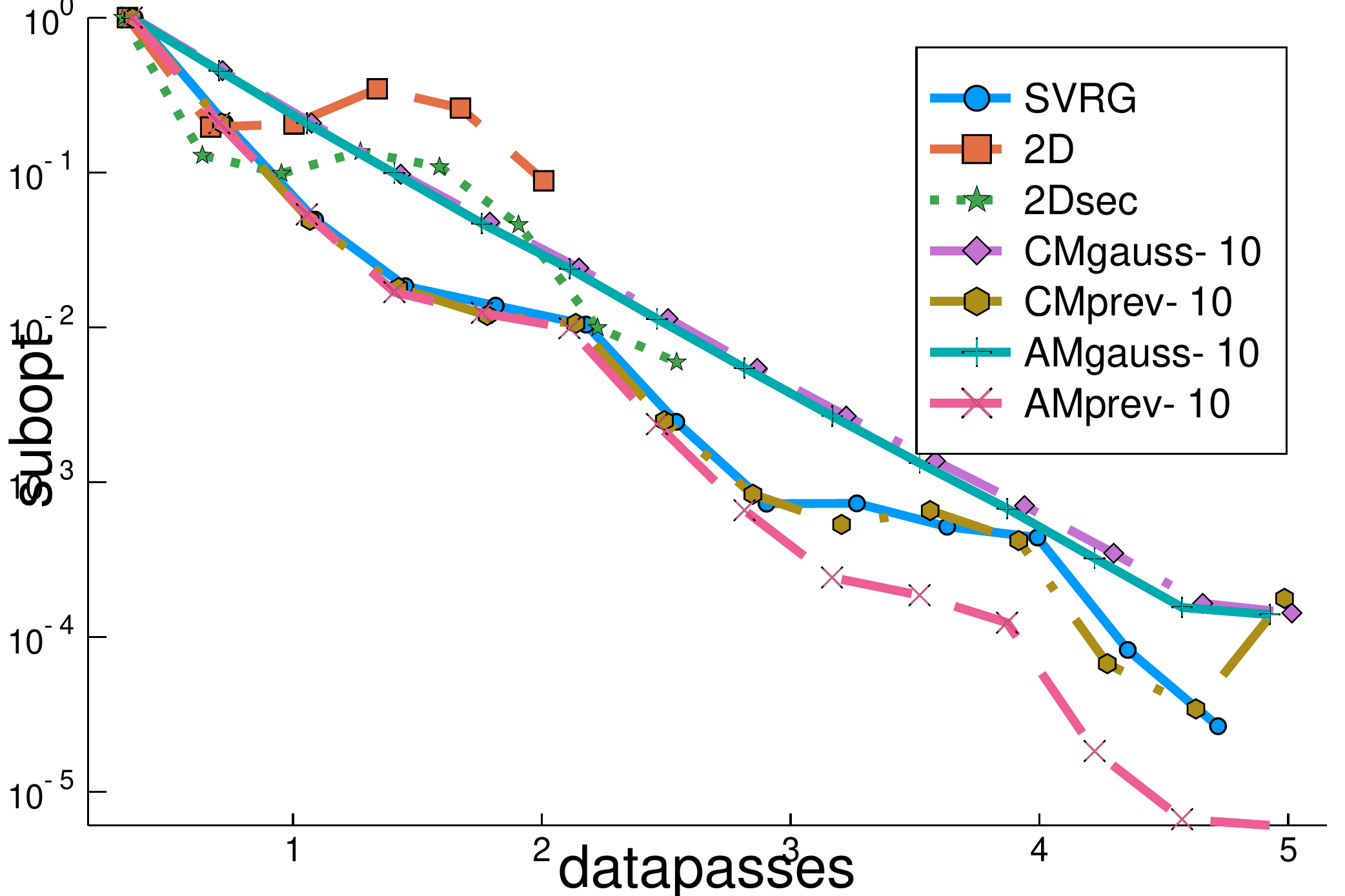}
		\caption{\texttt{rcv1-train} }
	\end{subfigure}
	\begin{subfigure}[t]{0.32\textwidth}
		\centering
		\includegraphics[width =  \textwidth]{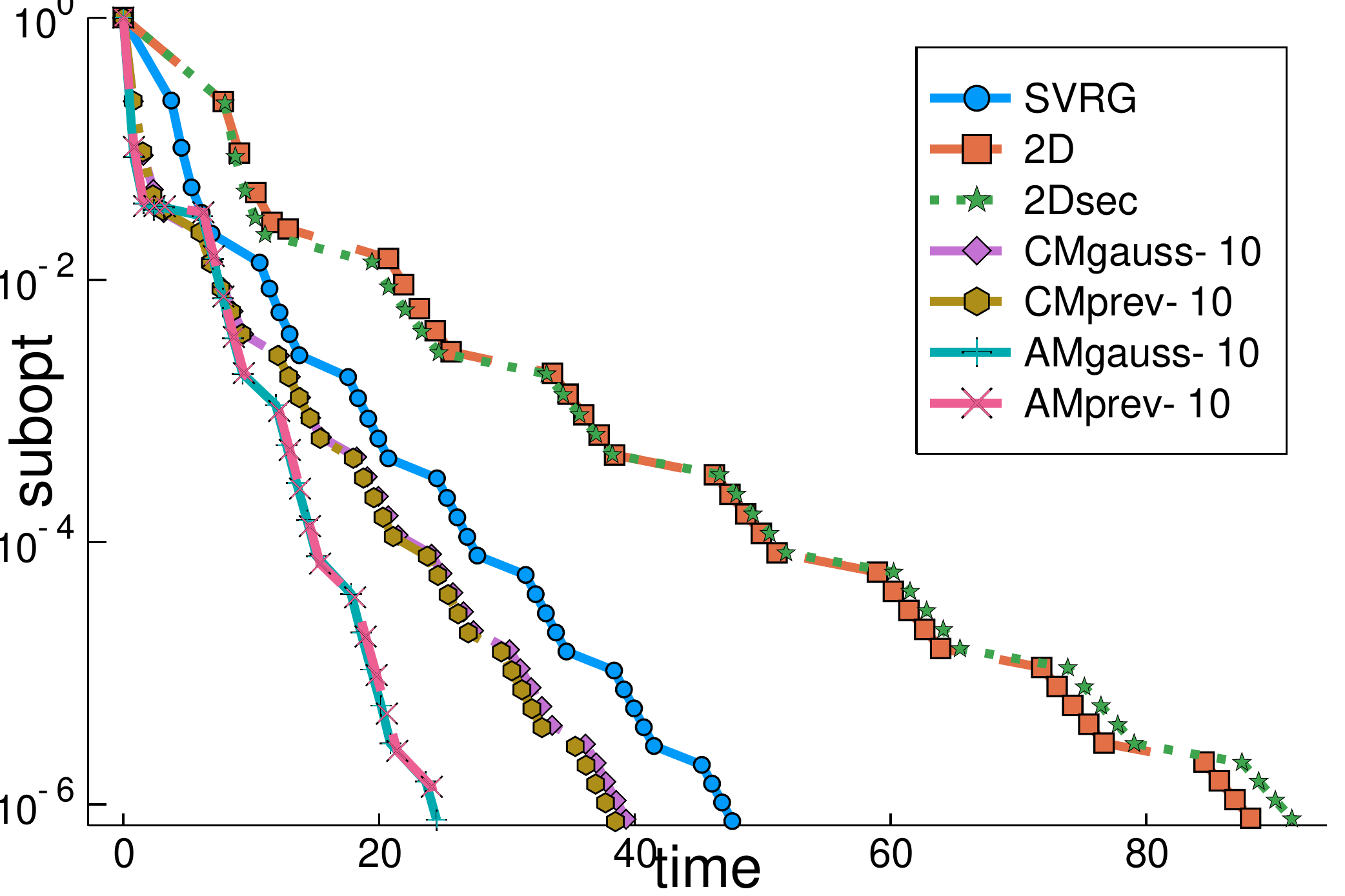} 
		\includegraphics[width =  \textwidth]{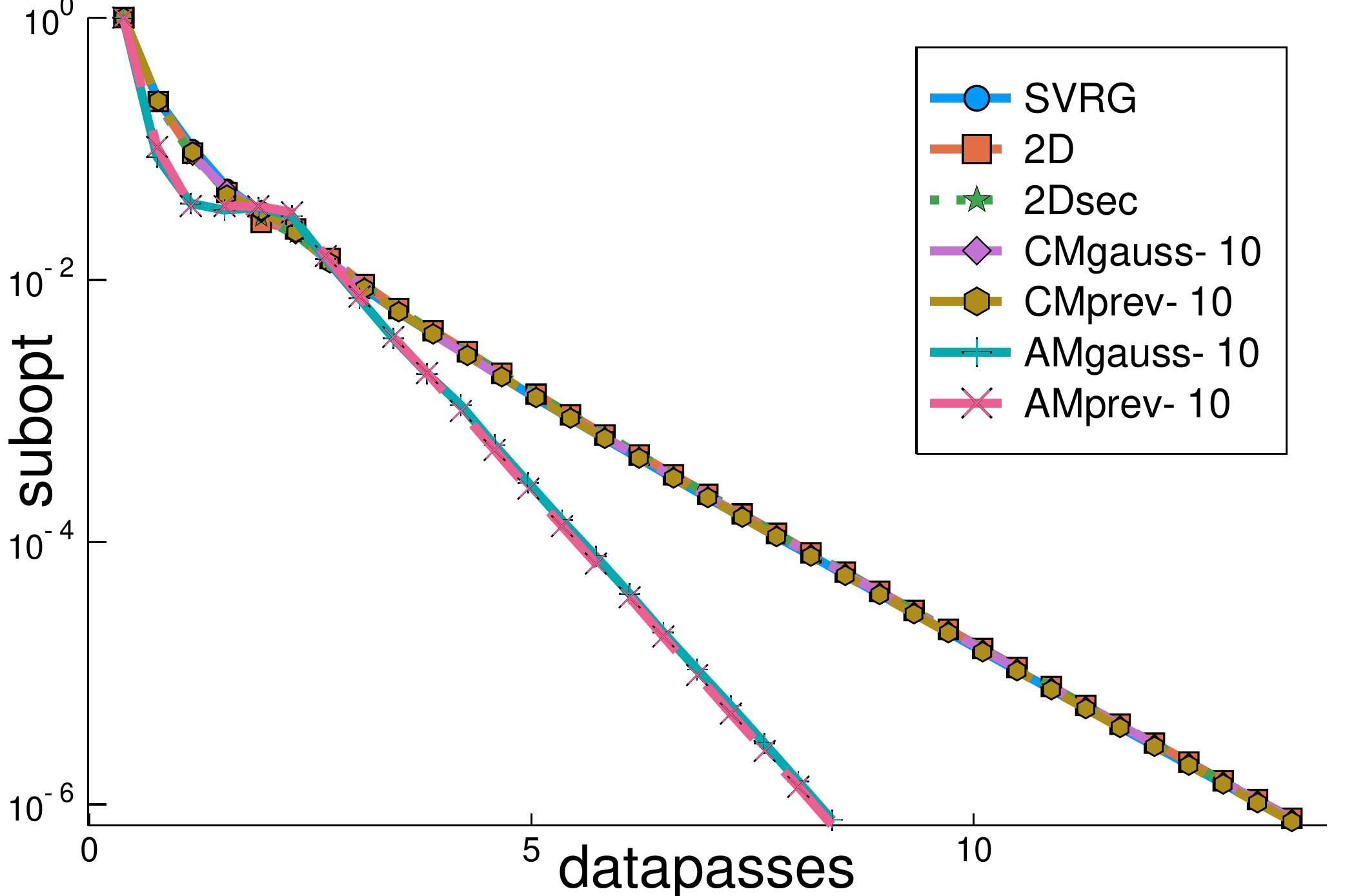}
		\caption{\texttt{gisette}}
	\end{subfigure}%
	\hspace{0.005\textwidth}
	\begin{subfigure}[t]{0.32\textwidth}
		\centering
		\includegraphics[width =  \textwidth]{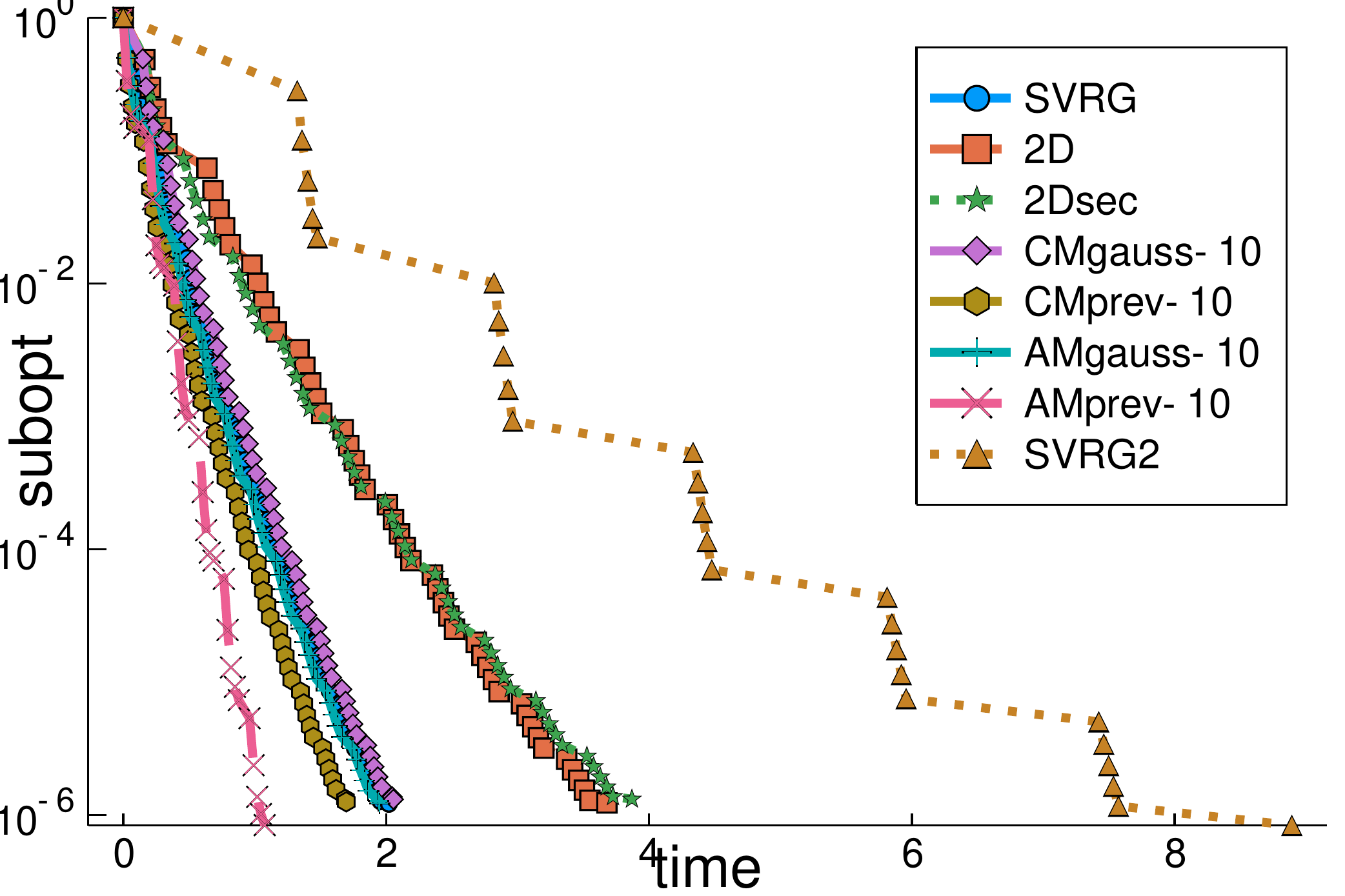} 
		\includegraphics[width =  \textwidth]{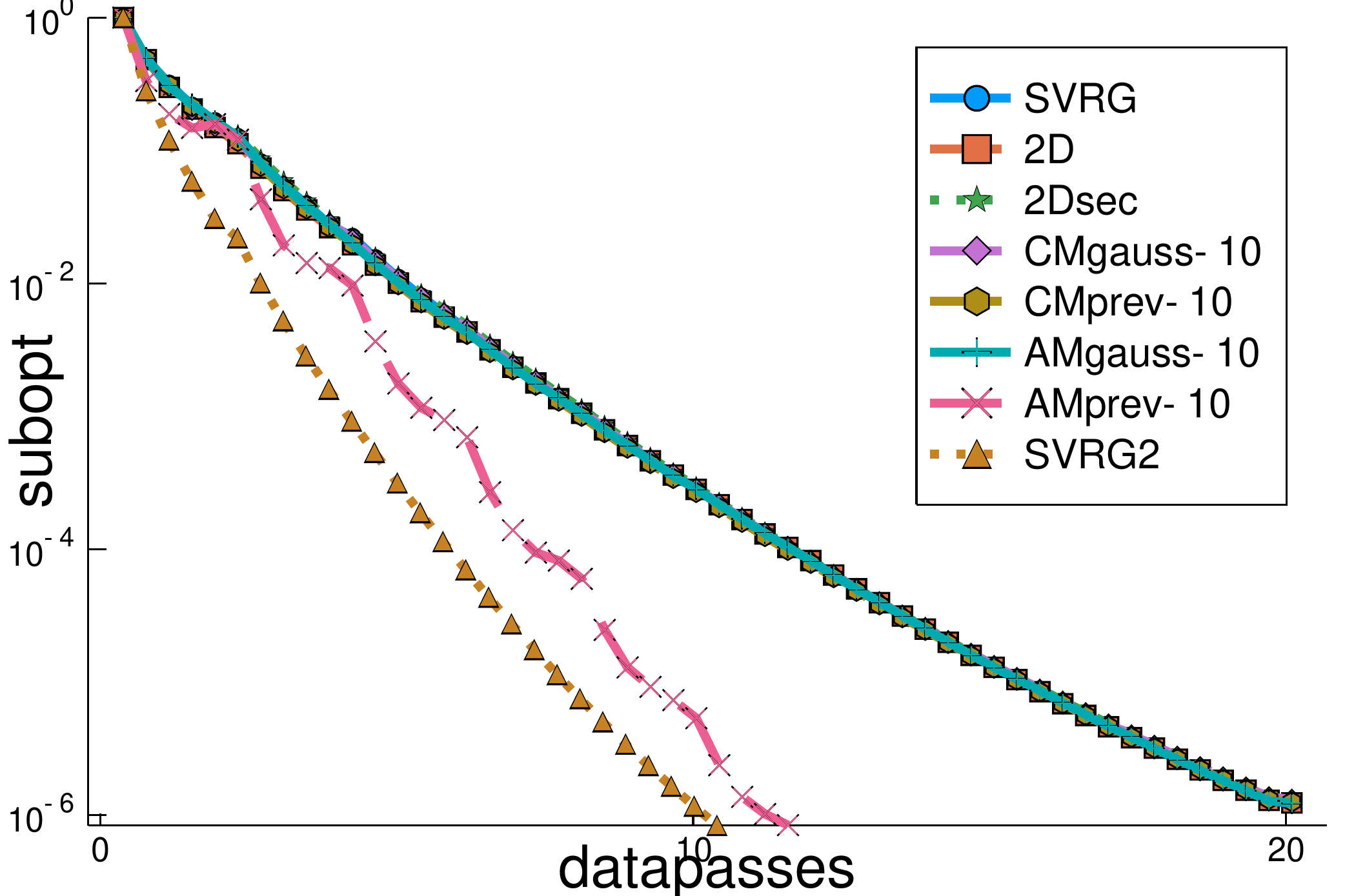}
		\caption{\texttt{madelon}}
	\end{subfigure}
	\caption{Performance of various SVRG-based methods on multiple LIBSVM problems:\\(a) \texttt{a9a} (b) \texttt{covtype} (c)  \texttt{mushrooms}  (d) \texttt{phishing} (e) \texttt{w8a}  (f) \texttt{rcv1-train}  (g) \texttt{gisette} (h) \texttt{madelon}.} \label{fig:LIBSVM}
\end{figure}

\textbf{Algorithms.}
We empirically tested the impact of using second-order information for the different Hessian approximations. All the code for our experiments was written in the julia progamming language and can be found at~\url{https://github.com/gowerrobert/StochOpt}. In particular, we tried \texttt{SVRG}, the original SVRG algorithm, \texttt{SVRG2}, which tracks the gradients using the full Hessian, \texttt{2D}, which tracks the gradients using the diagonal of the Hessian, \texttt{2Dsec}, which tracks the gradients using the robust secant equation~\eqref{eq:robsec}, \texttt{CM}, which tracks the gradients using the low-rank curvature matching approximation of the Hessian~\eqref{eq:embedleftright}, and \texttt{AM} which uses the low-rank action matching approximation of the Hessian~\eqref{eq:DFPinv} to track the gradients. We also tested two variants of each of the low rank approximations, namely \texttt{AMgauss, CMgauss} and \texttt{AMprev, CMprev}. For the \texttt{AMgauss} and \texttt{CMgauss} methods, the matrix $S \in \R^{n \times k}$ in~\eqref{eq:embedleftright} and~\eqref{eq:DFPinv} is a randomly generated Gaussian matrix. For the \texttt{AMprev} and \texttt{CMprev} methods, we construct $S$ using previous search directions according to~\eqref{eq:prevS}.

All the data we use are binary classification problems taken from LIBSVM, and we use a logistic loss with a squared-$\ell_2$ regularizer. We set the regularization parameter to $\lambda = \max_{i=1,\ldots, N} \norm{x_i}_2^2/4N$ in all experiments.  In each figure, the y-axis is the relative suboptimality gap given by $[f(\theta)-f(\theta^\ast)]/[f(\theta_0)-f(\theta^\ast)]$ on a log scale.

Convergence rate was observed both in terms of passes through the data and in terms of wall clock time, to account for the higher computational cost of our proposed methods.  Due to the high computational cost of \texttt{SVRG2}, we did not perform experiments on \texttt{SVRG2} for problems with dimension (number of features) over $5000$.
%

\textbf{Choice of the stepsize.}
Since convergence is guaranteed for batch gradient descent with a larger stepsize than for variance reducing stochastic methods, our hope is that using second-order information will allow the use of larger stepsizes while maintaining convergence. Thus, we ran a grid search over the stepsize for each method, trying values of $\left. 2^a \right/ L_{\max}$ for $a = 10,7, \dots, -7,-9,$ where $L_{\max} : = \max_{i=1,\ldots, N} \norm{x_i}_2^2 +\lambda.$
 
\begin{table}
\begin{tabular}{|l | cccccccc |} \hline
  	&	 \texttt{a9a} 	&		 \texttt{covtype} 	&		 \texttt{gisette} 	&		 \texttt{madelon} 	&		 \texttt{mushrooms} 	&		 \texttt{phishing}  	&		 \texttt{w8a}  	&		 \texttt{rcv1}		\\
\texttt{SVRG}   	&$	 2^ 8	$	&$	 2^ 9	$	&$	 2^ 8	$	&$	 2^ 7	$	&$	 2^ 8	$	&$	 2^ 6	$	&$	 2^ 8	$	&$	 2^ {10}	$	\\
\texttt{2D}     	&$	 2^ 8	$	&$	 2^ {10}	$	&$	 2^ 8	$	&$	 2^ 7	$	&$	 2^ 8	$	&$	 2^ 6	$	&$	 2^ 9	$	&$	 2^ {14}	$	        \\
\texttt{2Dsec}  	&$	 2^ 9	$	&$	 2^ {10}	$	&$	 2^ 8	$	&$	 2^ 7	$	&$	 2^ 8	$	&$	 2^ 6	$	&$	 2^ 9	$	&$	 2^ {13}	$	        \\
\texttt{CMgauss}        	&$	 2^ 8	$	&$	 2^ 9	$	&$	 2^ 8	$	&$	 2^ 7	$	&$	 2^ 8	$	&$	 2^ 6	$	&$	 2^ 9	$	&$	 2^ {9}	$	\\
\texttt{CMprev} 	&$	 2^ 8	$	&$	 2^ 8	$	&$	 2^ 8	$	&$	 2^ 7	$	&$	 2^ 7	$	&$	 2^ 6	$	&$	 2^ 9	$	&$	 2^ {10}	$	\\
\texttt{AMgauss}        	&$	 2^ 8	$	&$	 2^ 9	$	&$	 2^ 9	$	&$	 2^ 7	$	&$	 2^ 7	$	&$	 2^ 7	$	&$	 2^ 9	$	&$	 2^ {9}	$	\\
\texttt{AMprev} 	&$	 2^ 9	$	&$	 2^ 9	$	&$	 2^ 9	$	&$	 2^ 8	$	&$	 2^ 8	$	&$	 2^ 6	$	&$	 2^ 9	$	&$	 2^ {10}	$	\\
\texttt{SVRG2}  	&$	 2^{10}	$	&$	 2^ 9	$	&$	 2^ 9	$	&$	 2^ 8	$	&$	 2^ 9	$	&$	 2^ 6	$	&$	 2^ 9	$	&$	 2^ {10}	$	\\ \hline
\end{tabular}
\caption{Best empirical stepsize for each (problem, method) pair. In general, incorporating a Hessian approximation allows for larger stepsizes.}
\label{tab:stepsizes}
\end{table}

We expected our proposed covariates to provide better approximations to the true gradient and thus allow for greater stepsizes. This is the case in general, albeit not by a large amount, as can be seen in Table~\ref{tab:stepsizes}. Our proposed methods often allowed larger stepsizes, with the \texttt{SVRG2} method typically allowing the largest stepsize.

\begin{figure}
    \centering
    \begin{subfigure}[t]{0.235\textwidth}
        \centering
\includegraphics[width =  \textwidth]{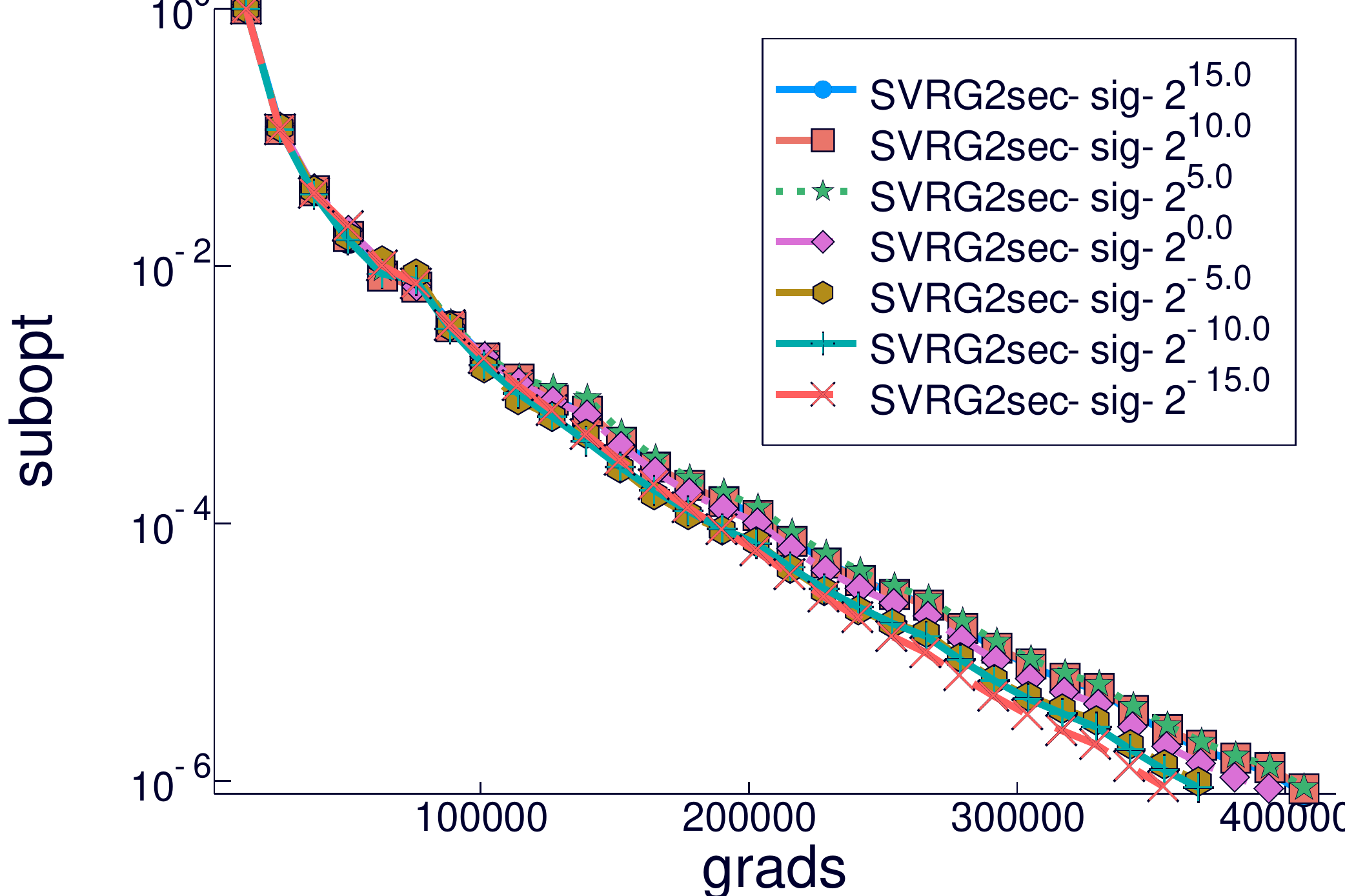}
        \caption{\texttt{a9a}}
    \end{subfigure}%
  \hspace{0.005\textwidth}
    \begin{subfigure}[t]{0.235\textwidth}
        \centering
\includegraphics[width =  \textwidth]{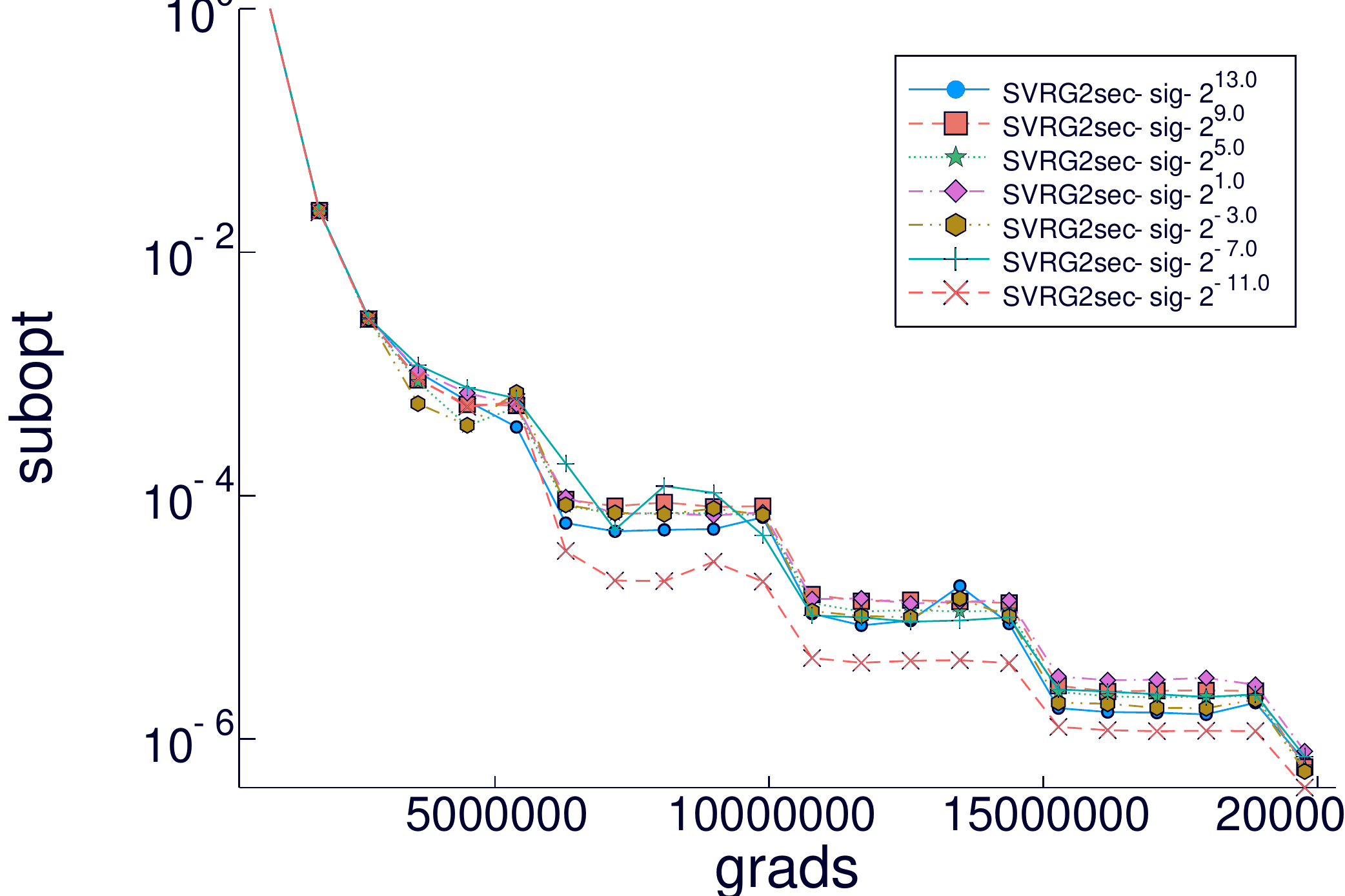}
        \caption{\texttt{covtype}}
    \end{subfigure}
      \hspace{0.005\textwidth}
    \begin{subfigure}[t]{0.235\textwidth}
        \centering
\includegraphics[width =  \textwidth]{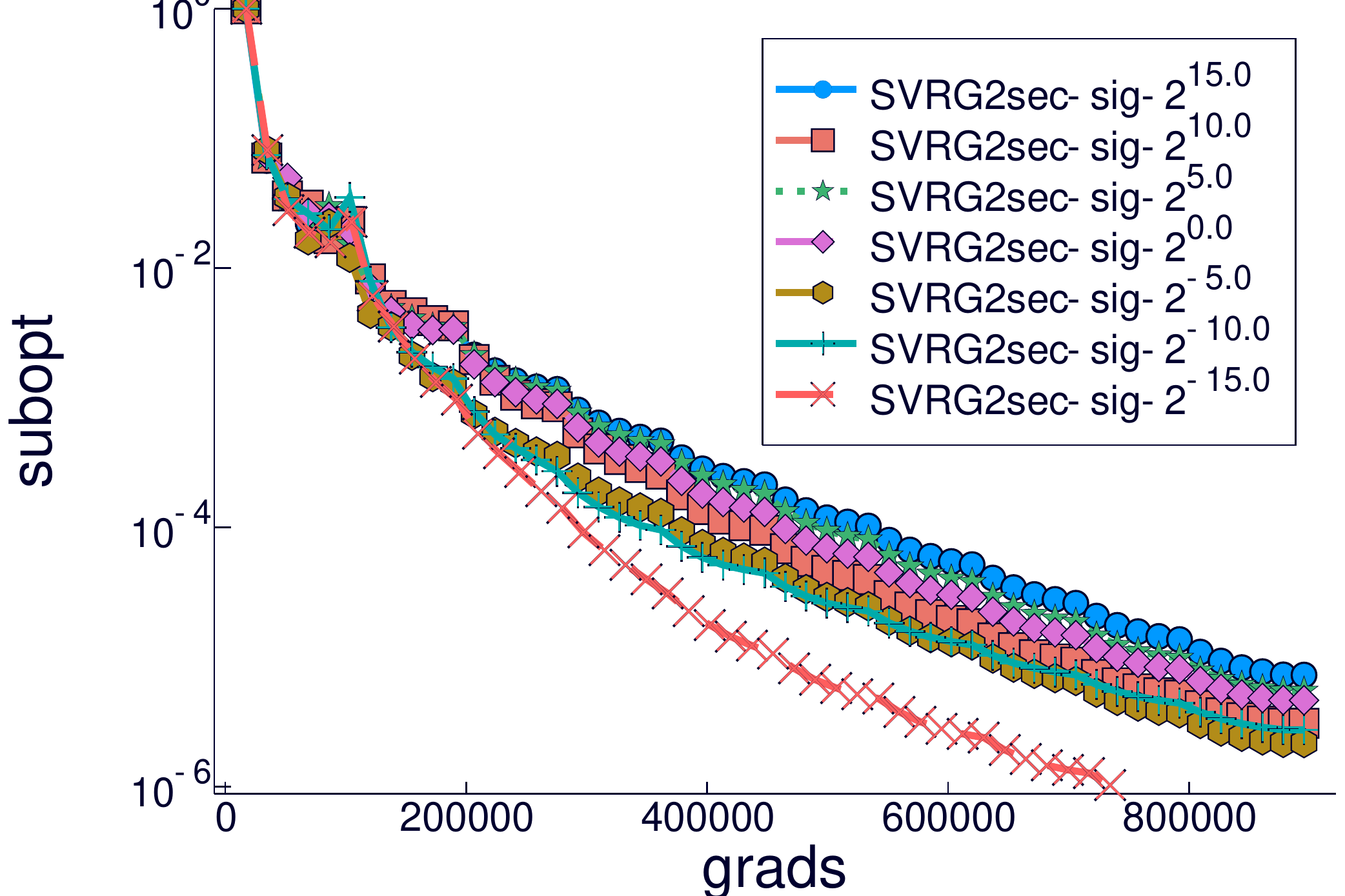}
        \caption{\texttt{phishing}}
    \end{subfigure}
      \hspace{0.005\textwidth}
    \begin{subfigure}[t]{0.235\textwidth}
        \centering
\includegraphics[width =  \textwidth]{mushrooms-sigmaexp}
        \caption{\texttt{mushrooms}}
    \end{subfigure}
    \caption{Performance of the SVRG2sec with different choices of $\sigma$ on: (a) \texttt{a9a} (b) \texttt{covtype} (c) \texttt{phishing} (d) \texttt{mushrooms}.} \label{fig:robust}
\end{figure}

\textbf{Results.}
Figure~\ref{fig:LIBSVM} shows the results on multiple datasets. We see that, in terms of passes through the data, all tracking methods consistently outperform \texttt{SVRG}, often by a large margin. In particular the \texttt{SVRG2} method often has the best performance in terms of datapasses with the exception of the \texttt{covtype} problem where the two methods based on diagonal approximations \texttt{2D} and \texttt{2Dsec} have the best performance. In terms of time taken, the \texttt{AMgauss} and \texttt{AMprev} are the most efficient methods, which indicates the superiority of the Hessian approximations based on action matching~\eqref{eq:actionmatch}.  For all these experiments, we chose a rank $k=10$.
\newpage
\textbf{Robustness of the diagonal approximation.}
Our robust secant equation has a hyperparameter, $\sigma^2$. Since the popularity of an optimization method depends as much of its ease of use as of its convergence rate, we tested the impact of $\sigma^2$ on the convergence speed. The results are presented in Fig.~\ref{fig:robust}. We can see that the impact is generally very limited and that our method is robust to the choice of $\sigma^2$. In all other experiments we set $\sigma^2 =0.1.$ 

\section{Conclusion}
Control variates are at the core of stochastic variance-reducing methods. We proposed the first dynamic method which updates these control variates at every timestep, leading to faster empirical and theoretical convergence as well as increased robustness. While we applied it to SVRG, we hope it can be extended to other methods such as SAGA.

We also showed that the Hessian could be used beyond reconditioning and that efficient approximations existed. Though reconditioning in stochastic variance-reducing methods has been explored, it remains to be seen whether using the Hessian for both reconditioning and gradient tracking simultaneously is beneficial.

Finally, we believe there is still more to do to bridge the gap between stochastic and batch methods and we hope this work will open the way to further analyses of the importance of control variates in optimization.

\subsection*{Acknowledgements}
The authors would like to acknowledge support from the Fondation de Sciences Math\'ematiques de Paris (FSMP) and from the European Research Council (ERC SEQUOIA).

\printbibliography

\appendix

\section{Theoretical Analysis: Proofs}\label{sec:theoretical-analysis:-proofs}

Following~\textcite{bubeck2015convex}, we consider a single epoch of SVRG and its extensions, that is, $\bar{\theta} \in \Theta$, and the iteration, started from $\theta_0 = \bar{\theta}$:
$$\theta_{t} = \Pi_\Theta \bigg( \theta_{t-1} - \gamma \Big[ f_{i_t}'(\theta_{t-1}) - z_{i_t}(\theta_{t-1}) + \frac{1}{N} \sum_{j=1}^N z_j(\theta_{t-1})\Big] \bigg),$$
with $i_t$ uniformly at random in $\{1,\dots,N\}$.

We also recall that, while the results are given using $R$ the radius of the data, they can be readily transposed to $L_{\max}$ using $L_{\max} = R^2$.

We have, with $\mathcal{F}_{t-1}$ representing the information up to time $t$:

\begin{eqnarray*}
\E { \| \theta_{t} - \theta_\ast\|^2 | \mathcal{F}_{t-1}}
& \leqslant & \E { \big\|  \theta_{t-1}- \theta_\ast - \gamma \Big[ f_{i_t}'(\theta_{t-1}) - z_{i_t}(\theta_{t-1}) + \frac{1}{N} \sum_{j=1}^N z_j(\theta_{t-1})\Big] \big\|^2 | \mathcal{F}_{t-1}}
\\
& & \hspace*{4cm} \mbox{ by contractivity of projections}, \\
& \leqslant &
\| \theta_{t-1} - \theta_\ast\|^2 - 2 \gamma F'(\theta_{t-1})^\top ( \theta_{t-1} - \theta_\ast) + \gamma^2 \| F'(\theta_{t-1})\|^2\\
& &
+ \frac{\gamma^2}{N} \sum_{i=1}^N
\Big\| f_{i}'(\theta_{t-1}) - z_{i}(\theta_{t-1})
- \frac{1}{N} \sum_{j=1}^N f_j'(\theta_{t-1}) + \frac{1}{N} \sum_{j=1}^N z_j(\theta_{t-1})
\Big\|^2
\\
& \leqslant & \| \theta_{t-1} - \theta_\ast\|^2 - 2 \gamma F'(\theta_{t-1})^\top ( \theta_{t-1} - \theta_\ast)   + \gamma^2 \| F'(\theta_{t-1})\|^2\\\
& &
+ \frac{\gamma^2}{N} \sum_{i=1}^N
\big\| f_{i}'(\theta_{t-1}) - z_{i}(\theta_{t-1})
\big\|^2, \mbox{ by bounding the variance by the second moment.}
\end{eqnarray*}

In the following sections, we provide proofs for several algorithms we consider in this paper.

\subsection{SVRG}

For regular SVRG (we provide the proof for completeness and because we need it later), we have: $z_i(\theta) = f_i'(\bar\theta)$ and we consider the bound
\begin{eqnarray*}
\frac{\gamma^2}{N} \sum_{i=1}^N
\big\| f_{i}'(\theta_{t-1}) - z_{i}(\theta_{t-1})
\big\|^2 & \leqslant    &
\frac{2\gamma^2}{N} \sum_{i=1}^N
\big\| f_{i}'(\theta_{t-1}) - f_{i}'(\theta_{\ast})
\big\|^2
+\frac{2\gamma^2}{N} \sum_{i=1}^N
\big\| f_{i}'(\bar \theta ) - f_{i}'(\theta_{\ast})
\big\|^2 \\
& \leqslant    &
 {2\gamma^2 R^2}  F'(\theta_{t-1})^\top ( \theta_{t-1} - \theta_\ast)
 +{2\gamma^2 R^2}  [ F(\bar \theta) - F(\theta_\ast)  ]
\end{eqnarray*}
leading to
\begin{eqnarray*}
\E { \| \theta_{t} - \theta_\ast\|^2 | \mathcal{F}_{t-1}}
& \leqslant &
\| \theta_{t-1} - \theta_\ast\|^2 - 2 \gamma F'(\theta_{t-1})^\top ( \theta_{t-1} - \theta_\ast)
+ \gamma^2 \| F'(\theta_{t-1})\|^2\\
& &
+ {2\gamma^2 R^2}  F'(\theta_{t-1})^\top ( \theta_{t-1} - \theta_\ast)
 +{2\gamma^2 R^2}  [ F(\bar \theta) - F(\theta_\ast)  ].
\end{eqnarray*}
Thus if $\gamma \leqslant 1 / ( 2R^2 + L )$, we get
\begin{eqnarray*}
\E { \| \theta_{t} - \theta_\ast\|^2 | \mathcal{F}_{t-1}}
& \leqslant &
\| \theta_{t-1} - \theta_\ast\|^2 -   \gamma[ F(\bar{\theta}_{t-1}) - F(\theta_\ast)  ] +{2\gamma^2 R^2}  [ F(\bar \theta) - F(\theta_\ast)  ].
\end{eqnarray*}
This implies that
\begin{eqnarray*}
\frac{1}{T}\sum_{t=1}^T
\E { F(\bar{\theta}_{t-1}) - F(\theta_\ast)  }
& \leqslant & \frac{1}{\gamma T} \|\bar\theta - \theta_\ast\|^2 +
{2\gamma R^2}  [ F(\bar \theta) - F(\theta_\ast)  ]
\\
\E { F \Big( \frac{1}{T} \sum_{t=1}^T \bar{\theta}_{t-1}\Big) - F(\theta_\ast)  }
& \leqslant & \Big( \frac{2}{\mu \gamma T}  +
{2\gamma  R^2} \Big) [ F(\bar \theta) - F(\theta_\ast)  ].
\end{eqnarray*}

This implies that if $\gamma = \frac{1}{4R^2}$ and $T \geqslant 8  / ( \gamma \mu) = \frac{32 R^2}{\mu}$, then
$$
\E { F \Big( \frac{1}{T} \sum_{t=1}^T \bar{\theta}_{t-1}\Big) - F(\theta_\ast)  }
 \leqslant  \frac{3}{4}  [ F(\bar \theta) - F(\theta_\ast)  ].
 $$
Thus, after $K = O ( \log \frac{1}{\varepsilon} ) $ epochs of SVRG we have attained the required precision, which makes an overall access to gradients
of $K N + K T = \big( N + \frac{R^2}{\mu}  \big)\log \frac{1}{\varepsilon}$.

\subsection{SVRG-2}
We assume that $\frac{4  \beta^2 R^ 4}{\alpha} D^2 \leqslant L$ and $\gamma = 1/ (4L)$.

In this situation, with no approximation, we have $z_i(\theta) = f_i'(\bar\theta) + f''_i(\bar\theta) ( \theta - \bar\theta)$ and:
\begin{eqnarray*}
 & & \frac{\gamma^2}{N} \sum_{i=1}^N
\big\| f_{i}'(\theta_{t-1}) - z_{i}(\theta_{t-1})
\big\|^2 \\
 & \leqslant    &
\frac{\gamma^2}{N} \sum_{i=1}^N R^2 \big[
\frac{\beta}{2} ( x_i^\top \theta_{t-1} - x_i^\top \bar\theta)^2
\big]^2  =
\frac{\gamma^2 \beta^2 R^ 2}{4 N} \sum_{i=1}^N   \
  ( x_i^\top ( \theta_{t-1} -  \bar\theta) )^4 \mbox{ using the bound on } \varphi''',\\
& \leqslant    &
\frac{\gamma^2 \beta^2 R^ 2}{N} \sum_{i=1}^N
\big[ 2   ( x_i^\top ( \theta_{t-1} -  \theta_\ast) )^4
+ 2   ( x_i^\top ( \theta_{\ast} -  \bar\theta) )^4 \big]
 \\
& \leqslant    &
\frac{\gamma^2 \beta^2 R^ 2}{N} \sum_{i=1}^N
\big[ 2R^2 \| \theta_{t-1} -  \theta_\ast \|^2  ( x_i^\top ( \theta_{t-1} -  \theta_\ast) )^2
+ 2 R^2 \| \bar\theta -  \theta_\ast \|^2  ( x_i^\top ( \bar\theta -  \theta_\ast) )^2\big]
\mbox{ using } \|x_i \| \leqslant R, \\
& \leqslant    &
\frac{2 \gamma^2 \beta^2 R^ 4}{N} \| \theta_{t-1} -  \theta_\ast \|^2 \sum_{i=1}^N
  ( x_i^\top ( \theta_{t-1} -  \theta_\ast) )^2
  +\frac{2 \gamma^2 \beta^2 R^ 4}{N} \|  \bar\theta -  \theta_\ast \|^2 \sum_{i=1}^N
  ( x_i^\top (  \bar\theta  -  \theta_\ast) )^2
\\
& \leqslant    &
\frac{4 \gamma^2 \beta^2 R^ 4}{\alpha } \| \theta_{t-1} -  \theta_\ast \|^2
[ F(\theta_{t-1}) - F(\theta_\ast) ]
+
\frac{4 \gamma^2 \beta^2 R^ 4}{\alpha } \| \bar\theta  -  \theta_\ast \|^2
[ F( \bar\theta) - F(\theta_\ast) ]
\mbox{ using } \varphi'' \geqslant \alpha,\\
& \leqslant    &
\frac{4 \gamma^2 \beta^2 R^ 4}{\alpha } D^2
[ F(\theta_{t-1}) - F(\theta_\ast) ]
+
\frac{4 \gamma^2 \beta^2 R^ 4}{\alpha  } D^2
[ F( \bar\theta) - F(\theta_\ast) ], \mbox{ using the compactness of } \Theta.
\end{eqnarray*}

With our assumptions, we have  $ { \gamma \bigg( L + \frac{4  \beta^2 R^ 4}{\alpha} D^2
\bigg) \leqslant 1} $, and we get that
\begin{eqnarray*}
\E { \| \theta_{t} - \theta_\ast\|^2 | \mathcal{F}_{t-1}}
& \leqslant &
\| \theta_{t-1} - \theta_\ast\|^2 - \gamma F'(\theta_{t-1})^\top ( \theta_{t-1} - \theta_\ast)
+ \frac{4 \gamma^2 \beta^2 R^ 4}{\alpha  } D^2    [ F(\bar \theta) - F(\theta_\ast)  ] \\
& \leqslant &
\| \theta_{t-1} - \theta_\ast\|^2 - \gamma  [ F(\theta_{t-1}) - F(\theta_\ast)  ]
+ \frac{4 \gamma^2 \beta^2 R^ 4}{\alpha  } D^2    [ F(\bar \theta) - F(\theta_\ast)  ].
\end{eqnarray*}

This leads to, with  $\boxed{T \geqslant 4  / ( \mu \gamma ) = \frac{ 16 L}{\mu} }$
and using ${ \gamma \bigg( \frac{4  \beta^2 R^ 4}{\alpha} D^2
\bigg) \leqslant 1/2} $,
\begin{eqnarray*}
\E { F \Big( \frac{1}{T} \sum_{t=1}^T \bar{\theta}_{t-1}\Big) - F(\theta_\ast)  }
& \leqslant &
\bigg(
\frac{2}{\mu \gamma T} + \frac{4 \gamma \beta^2 R^ 4}{\alpha  } D^2
\bigg)
[ F(\bar \theta) - F(\theta_\ast)  ]
\\
 & \leqslant  &  \frac{3}{4}  [ F(\bar \theta) - F(\theta_\ast)  ].
 \end{eqnarray*}
Thus, after $K = O ( \log \frac{1}{\varepsilon} )$ epochs of SVRG we have attained the required precision, which makes an overall access to gradients
of $K N + K T = \big( N + \frac{L}{\mu}   \big)\log \frac{1}{\varepsilon}$.

\subsection{Stability of SVRG-2}
If we make no compactness assumption on $\Theta$, then we have:
\begin{eqnarray*}
 & & \frac{\gamma^2}{N} \sum_{i=1}^N
\big\| f_{i}'(\theta_{t-1}) - z_{i}(\theta_{t-1})
\big\|^2 \\
& \leqslant &
\frac{2 \gamma^2}{N} \sum_{i=1}^N
\big\| f_{i}'(\theta_{t-1}) - f_{i}(\bar{\theta})
\big\|^2
+ \frac{2 \gamma^2}{N} \sum_{i=1}^N
\big\| f_{i'}''(\bar\theta)(\theta_{t-1} - \bar{\theta} )
\big\|^2 \\
& \leqslant &  {2\gamma^2 R^2}  F'(\theta_{t-1})^\top ( \theta_{t-1} - \theta_\ast)
 +{2\gamma^2 R^2}  [ F(\bar \theta) - F(\theta_\ast)  ]
\mbox{ from the SVRG proof },
\\
& & + \frac{2 \gamma^2}{N} \sum_{i=1}^N R^2
\big\| x_i^\top (\theta_{t-1} - \bar{\theta} )
\big\|^2
\\
& \leqslant &  {2\gamma^2 R^2}  F'(\theta_{t-1})^\top ( \theta_{t-1} - \theta_\ast)
 +{2\gamma^2 R^2}  [ F(\bar \theta) - F(\theta_\ast)  ]
\\
& & \frac{2\gamma^2 R^2}{\alpha}  F'(\theta_{t-1})^\top ( \theta_{t-1} - \theta_\ast)
 + \frac{2\gamma^2 R^2}{\alpha}  [ F(\bar \theta) - F(\theta_\ast)  ]
 \\
 & \leqslant & \frac{4\gamma^2 R^2}{\alpha}  F'(\theta_{t-1})^\top ( \theta_{t-1} - \theta_\ast)
 + \frac{4\gamma^2 R^2}{\alpha}  [ F(\bar \theta) - F(\theta_\ast)  ]
 \end{eqnarray*}

 Thus, if we take the smaller step-size $\gamma = \frac{\alpha}{8 R^2}$ and $T = \frac{64 R^2}{\alpha \mu}$, we get the same convergence.
\end{document}